\documentclass[12pt]{article}

\usepackage{latexsym}
\usepackage{calc}
\usepackage{datetime}
\usepackage{amssymb, amsmath, amsfonts, listings, mathrsfs} 
\usepackage{xcolor} 
\usepackage[normalem]{ulem} 

\usepackage{graphicx}
\usepackage[labelfont=bf]{caption}

\newcounter{hours}\newcounter{minutes}

\newcommand{\stkout}[1]{\ifmmode\text{\sout{\ensuremath{#1}}}\else\sout{#1}\fi} 

\def\nr{\par \noindent}
\def\diag{{\rm diag \,}}
\def\Diag{{\rm Diag \,}}
\def\Def{\stackrel{\mathrm{def}}{=}}

\def\inter{{\rm int \,}}

\def\dom{{\rm dom \,}}
\def\vf{\varphi}
\def\vk{\varkappa}
\def\beq{\begin{equation}}
\def\eeq{\end{equation}}
\newcommand{\rint}{{\rm rint\,}}

\def\C{\mathbb{C}}
\def\R{\mathbb{R}}
\def\E{\mathbb{E}}

\def\S{\mathbb{S}}

\def\BI{\begin{itemize}}
\def\EI{\end{itemize}}

\newcommand{\SetEQ}{\setcounter{equation}{0}}
\newcommand{\refLE}[1]{\ensuremath{\stackrel{(\ref{#1})}{\leq}}}
\newcommand{\refLEQ}[2]{\ensuremath{\stackrel{(\ref{#1}),(\ref{#2})}{\leq}}}

\newcommand{\refLEI}[2]{\ensuremath{\stackrel{(\ref{#1}_{#2})}{\leq}}}
\newcommand{\refEQ}[1]{\ensuremath{\stackrel{(\ref{#1})}{=}}}
\newcommand{\refEQI}[2]{\ensuremath{\stackrel{(\ref{#1}_{#2})}{=}}}
\newcommand{\refEQQ}[2]{\ensuremath{\stackrel{(\ref{#1}),(\ref{#2})}{=}}}
\newcommand{\refGE}[1]{\ensuremath{\stackrel{(\ref{#1})}{\geq}}}

\newtheorem{theorem}{Theorem}
\newtheorem{lemma}{Lemma}
\newtheorem{corollary}{Corollary}

\newtheorem{assumption}{Assumption}
\newtheorem{definition}{Definition}

\newtheorem{example}{Example}
\newtheorem{remark}{Remark}
\newcommand{\proof}{\bf Proof: \rm \nr}
\newcommand{\qed}{\hfill $\Box$ \nr \medskip}
\newcommand{\half}{\mbox{${1 \over 2}$}}

\def\ba{\begin{array}}
\def\ea{\end{array}}
\def\beann{\begin{eqnarray*}}
\def\eeann{\end{eqnarray*}}
\def\bea{\begin{eqnarray}}
\def\eea{\end{eqnarray}}

\def\BT{\begin{theorem}}
\def\ET{\end{theorem}}
\def\BL{\begin{lemma}}
\def\EL{\end{lemma}}
\def\BC{\begin{corollary}}
\def\EC{\end{corollary}}
\def\BE{\begin{example}}
\def\EE{\end{example}}
\def\BD{\begin{definition}}
\def\ED{\end{definition}}
\def\BR{\begin{remark}}
\def\ER{\end{remark}}
\def\BAS{\begin{assumption}}
\def\EAS{\end{assumption}}
\def\BI{\begin{itemize}}
\def\EI{\end{itemize}}

\def\BMP{\begin{minipage}{9.5cm}}
\def\EMP{\end{minipage}}
\def\MPT{\begin{minipage}{11.5cm}}
\def\EPT{\end{minipage}}

\def\la{\langle}
\def\ra{\rangle}

\def\QF{\hspace{5ex} \Box}
\def\QR{\hfill \Box}

\title{
Asymmetric Long-Step Primal-Dual Interior-Point \\
Methods with Dual Centering 
}

\author{Yurii Nesterov
\thanks{Corvinus Centre for Operations Research, Corvinus Institute for Advanced Studies, Corvinus University of Budapest, and School of Data Sciences (Chinese University of Hong Kong (Shenzhen); part time). \newline Professor emeritus at CORE/INMA, UClouvain, Belgium. Email: Yurii.Nesterov@uclouvain.be. \newline
This research was supported by the National Research, Development and Innovation Office (NKFIH) under grant number 2024-1.2.3-HU-RIZONT-2024-00030.}}

\date{\normalsize 
February 16, 2025
}

\begin{document}
\maketitle

\abstract{In this paper, we develop a new asymmetric framework for solving primal-dual problems of Conic Optimization by Interior-Point Methods (IPMs). It allows development of efficient methods for problems, where the dual formulation is simpler than the primal one. The problems of this type arise, in particular, in Semidefinite Optimization (SDO), for which we propose a new method with very attractive computational cost. Our long-step predictor-corrector scheme is based on centering in the dual space. It computes the affine-scaling  predicting direction by the use of the dual barrier function, controlling the tangent step size by a functional proximity measure. We show that for symmetric cones, the search procedure at the predictor step is very cheap.

In general, we do not need sophisticated Linear Algebra, restricting ourselves only by Cholesky factorization. However, our complexity bounds correspond to the best known polynomial-time results. Moreover, for symmetric cones the bounds automatically depend on the minimal barrier parameter between the primal or the dual feasible sets. We show by SDO-examples that the corresponding gain can be very big.

We argue that the dual framework is more suitable for adjustment to the actual complexity of the problem. As an example, we discuss some classes of SDO-problems, where the number of iterations is proportional to the square root of the number of linear equality constraints. Moreover, the computational cost of one iteration there is similar to that one for Linear Optimization. We support our theoretical developments by  preliminary but encouraging numerical results with randomly generated SDO-problems of different size.}

\vspace{2ex}\noindent
{\bf Keywords:} Conic optimization, self-concordant barriers, polynomial-time interior-point methods, primal-dual methods, predictor-corrector methods, Semidefinite Optimization.

\vspace{2ex}\noindent
{\bf Mathematics Subject Classification:} 90C51, 90C25, 90C99



\section{Introduction}
\setcounter{equation}{0}

\vspace{1ex}\noindent
{\bf Motivation.} After the seminal papers \cite{Kar,Jim,Clovis}, the end of 20th century was marked by an enormous research activity related to polynomial-time Interior-Point Methods. Thousands of papers were published in a short period of time, resulting in a completely new framework for development of  efficient optimization schemes. The first really successful algorithms were proposed for Lineaar Optimization in the form of primal-dual predictor-corrector methods (see \cite{Mer,MTY,ZZ}, and \cite{YZ} for the recent developments). Later, they were extended onto important family of symmetric cones \cite{Ali,Self1,Self2,Todd,Boyd} and onto the general cones in the framework of self-concordant functions (see \cite{NN} and \cite{Mosek} for the recent results). The interested reader can find different expositions of this theory in many monographs and surveys \cite{BN,PS,RB,RVT,MWright,Wright,Ye}.

Today, the primal-dual formulation of conic problems \cite{NN,BN} remains one of the most powerful framework in Nonlinear Convex Optimization. Its high performance is usually explained by a full symmetry between the primal and the dual problems. Indeed, these formulations help each other in constructing powerful primal-dual methods by getting an additional information from a conjugate pair of self-concordant barriers for the primal and the dual cones. The most efficient algorithms are predictor-corrector methods \cite{Mer,MTY,LongStep}, which follow approximately the primal-dual central path by alternating the predictor step with several correcting iterations.

However, the superiority of these methods is not for free. For practical applications, we need to assume computability both for the primal and the dual barrier functions. For general cones, this ability is quite rare. Indeed, in order to compute the value of the dual barrier function, we need to find in a closed form a minimum of the sum of the primal barrier and a linear function with arbitrary coefficients. For nontrivial barrier functions, this is usually impossible. Of course, we could try to evaluate the dual barrier numerically. However, any inaccuracy in the computation of its value and/or its derivatives is damaging for the optimization schemes.

Thus, at this moment the classical primal-dual methods work well only for the self-dual barriers, which correspond mainly to symmetric cones.\footnote{For practical applications, the most important cones are the positive orthant, the cone of positive semidefinite matrices and the Lorentz cone.}

The above asymmetry in the computational cost of the primal and dual barrier functions served as the main motivation for the paper \cite{TW}, where the role of the dual barrier function was significantly reduced. It was suggested to perform the correction process only by the primal barrier function. This function and its derivatives are also used for computing the predicting direction. The dual barrier is necessary only for evaluating the {\em functional proximity measure} \cite{LongStep}, which can accept approximate values of the dual barrier.

However, it appears that the computational difficulties with the dual barrier are not the only source of asymmetry in the primal-dual formulation. Even for symmetric cones, the complexity of primal and dual problems can be very different. In Semidefinite Programming this can be seen even at the structural level. Indeed, the primal variables there are symmetric $n \times n$-matrices, and the dual variables are vectors from $\R^m$, where $m$ is the number of equality constraints. Thus the size of the primal space is much bigger than that of the dual space. Moreover, it appears that in some cases the complexity of the dual feasible set, measured by the value of barrier parameter, is much smaller than the complexity of the primal feasible set. In Section \ref{sc-Easy}, we give examples of SDO-problems where the dual barrier parameter is equal to $m$, which can be much smaller than the dimension of matrices $n$.

The main goal of this paper is the development of a primal-dual framework for solving problems of Conic Optimization, based mainly on the {\em dual} barrier function. In fact, in many applications the dual formulation is simpler. To give a motivating example, we prove that just a standard assumption on boundedness of the dual feasible set ensures the value of dual barrier parameter to be smaller than that one for the primal feasible set.

Hence, we present a dual counterpart of the primal technique developed in \cite{TW}. One of the main applications of our developments is Semidefinite Optimization, where we suggest a new way for generating the prediction directions. It is much less expensive than the standard one based on computing the scaling point for a given feasible primal-dual pair \cite{Self1,Self2}. In particular, in our approach, we do not need to compute the square roots of symmetric matrices. We use only the straightforward Cholesky factorizations.

\vspace{1ex}\noindent
{\bf Contents.} In Section \ref{sc-Not}, we introduce the main notation and recall the basic facts from the theory of self-concordant functions. In Section \ref{sc-PDPaths}, we define the primal and dual central paths and study their relations with the functional proximity measure \cite{LongStep}. Also, we discuss the importance of {\em scaling point}, whose Heasian transforms the dual feasible point into the primal one. Using this Hessian, we define the primal-dual {\em Affine Scaling Direction} for the predictor step. An important identity (\ref{eq-LNorms}) for the local norms of the restricted barriers for the primal and the dual feasible sets is one of our main arguments supporting the importance of asymmetric treatment of the primal and the dual problems. The section is concluded with an example of SDO-problem with maximal primal-dual asymmetry.

In Section \ref{sc-BDS}, we present the main element of our construction, the {\em Dual Gambit Step}. It is defined by the point obtained in the dual correction process. This point serves as a {\em scaling point}, connecting some feasible primal and dual variables, defined explicitly by a simple rule. These points remain in a close neighborhood of the primal-dual central path.

In Section \ref{sc-ASD}, we define the {\em Affine-Scaling Direction} and show that the growth of the functional proximity measure along this direction is very slow. This direction plays a key role in the Dual Primal-Dual Predictor-Corrector Method, explained in Section \ref{sc-DBMet}. For general cones, its complexity bound is of the order $O(\sqrt{\nu} \ln {1 \over \epsilon})$, where $\nu$ is the parameter of the barrier for the primal feasible cone and $\epsilon$ is the desired accuracy of the solution.

In Section \ref{sc-Sym}, we present a variant of method of Section \ref{sc-DBMet} for symmetric cones. For these cones, it is possible to represent the primal-dual functional proximity measure completely in terms of the dual barrier function. This representation leads to much sharper complexity bounds. In particular, we show that the total number of predictor steps in the method is defined by the value of barrier parameter for the dual feasible set. As compared with another attempt \cite{NTun} of solving the dual problem by a long-step path-following scheme, we use now much simpler proximity measure.

In Section \ref{sc-Easy} we analyse performance of our method at the problems of Semidefinite Optimization. We show that for some problems, the complexity of the dual feasible set is much lower than that of the primal set. In particular, the value of parameter of the dual barrier function can be bounded by the sum of the ranks of the equality constraints. It can be much smaller than the dimension of matrices. It is important that our method adjusts automatically to the correct value of barrier parameter. We justify also a short representation of the dual barrier function, which can be useful in some problem settings.

In Section \ref{sc-Num}, we present preliminary computational results for problems of Low-Rank Quadratic Interpolation, which have complexity bound of $O(\sqrt{m} \ln {1 \over \epsilon})$ iterations. Section~\ref{sc-Conc} contains Conclusion. And in Section \ref{sc-App} (Appendix) we show how to bound the growth of functional proximity measure by the barrier parameter of primal feasible set.

In the sequel, the references onto specific parts of the displayed relations are indicated by the corresponding index (e.g. (\ref{eq-FDual}$_3$)).

\section{Notations and Generalities}\label{sc-Not}
\SetEQ

Let $\E$ be a finite-dimensional space and $\E^*$ be the space of linear functions on $\E$. For $x \in \E$ and $s \in \E^*$, denote by $\la s, x \ra$ the value of function $s$ at $x$. We use the same notation $\la \cdot, \cdot \ra$ for different spaces. Thus, its actual meaning is defined by the context.
For a linear operator $A: \E_1 \to \E^*_2$, we denote by $A^*$ its {\em adjoint operator}:
\beq\label{def-AD}
\ba{rcl}
\la A x, y \ra & = & \la A^* y, x \ra, \quad x \in \E_1, \; y \in  \E_2.
\ea
\eeq
Thus, $A^*: \E_2 \to \E_1^*$. The operator $A: \E \to \E^*$ is positive semidefinite (notation $A \succeq 0$) if $\la A x, x \ra \geq 0$ for all $x \in \E$. We write $A \succeq B$ if $A-B \succeq 0$. 

In case of $\E = \R^n$, we have $\E^* = \R^n$ and
$$
\ba{rcl}
\la s, x \ra & = & \sum\limits_{i=1}^n s^{(i)} x^{(i)}, \quad \| x \| \; \Def \; \la x, x \ra^{1/2}, \quad s, x \in \R^n.
\ea
$$
Vector $e \in \R^n$ is the vector of all ones, and $e_i$, $i = 1, \dots , n$, are coordinate vectors in $\R^n$.

For vectors $x, s \in \R^n$, all arithmetic operations like $x \pm s$ or $x s$ are understood in the coordinate-wise sense. The same is true for the partial orderings.

Denote by $\R^{n \times m}$ the space of $n \times m$-matrices with the Frobenius scalar product 
$$
\ba{rcl}
\la S, X \ra & = & \sum\limits_{i=1}^n \sum\limits_{j=1}^m S^{(i,j)} X^{(i,j)}, \quad S, X \in \R^{n \times m}.
\ea
$$
For $x \in \R^n$, notation $X \equiv D(x) \in \R^{n \times n}$ is used for the diagonal matrix with $x$ at the main diagonal. Thus, the identity matrix is $I_n \equiv D(e)$.

Notation ${\cal C}^k(\E)$, $k \geq 1$ is used for the class of closed functions on $\E$ with open domain $\dom f$, which admit $k$th directional derivative  $D^kf(x)[h]^k$ along any direction $h \in \E$ at any point $x \in \dom f$. Function $f \in {\cal C}(\E)$ is differentiable on $\dom f$ (notation $f \in \C(\E)$) if
$$
\ba{rcl}
DF(x)[h] & = & \la \nabla f(x), h \ra, \quad x \in \dom f, \; h \in \E.
\ea
$$
Thus, the gradient $\nabla f(x)$ belongs to $\E^*$, $x \in \dom f$.
Function $f \in {\cal C}^2(\E)$ is twice differentiable on $\dom f$ (notation $f \in \C^2(\E)$) if
$$
\ba{rcl}
D^2F(x)[h]^2 & = & \la \nabla^2 f(x)h, h \ra, \quad x \in \dom f, \; h \in \E.
\ea
$$
Thus, for any  $x \in \dom f$, the linear operator (Hessian) $\nabla^2 f(x)$ maps $\E$ to $\E^*$.

Let us recall some facts from the theory of self-concordant functions (e.g. \cite{Lect,NN}). 

\BD 
Function $f \in {\cal C}^3(\E)$ is called {\em self-concordant} (SCF) if its third directional derivative is bounded by an appropriate power of the second derivative:
\beq\label{def-SCF}
\ba{rcl}
D^3f(x)[h]^3 & \leq & 2 \Big( D^2 f(x)[h]^2 \Big)^{3/2}, \quad x \in \dom f, \; h \in \E.
\ea
\eeq
\ED
If $\dom f$ contains no straight line, then its Hessian is positive definite at any $x \in \dom f$. In this paper, we mainly deal with such functions.

For any SCF, its dual function $f_*(s) = \max\limits_{x \in \dom f} [ - \la s, x \ra - f(x)]$ is also self-concordant. Note that $- \nabla f(x) \in \dom f_*$ for all $x \in \dom f$, and
$$
\ba{rcl}
\nabla f_*(s) & = & - x(s) \; \Def\; \arg\min\limits_{x \in \dom f} [ - \la s, x \ra - f(x) ], \quad s \in \dom f_*
\ea
$$
For any $x \in \dom f$ and $s \in \dom f_*$, we have
\beq\label{eq-FDual}
\ba{rclrcl}
\nabla f(-\nabla f_*(s)) & = & s, & \nabla^2 f(-\nabla f_*(s)) & = & [ \nabla^2 f_*(s)]^{-1},\\
\\
\nabla f_*(-\nabla f(x)) & = & x &
\nabla^2 f_*(-\nabla f(x)) & = & [ \nabla^2 f(x)]^{-1}.
\ea
\eeq

One of the main properties of SCF is that the {\em Dikin Ellipsoid} defined as
$$
\ba{rcl}
W^f_r(x) & = & \{ u \in \E:\; \la \nabla^2 f(x)(y-x),y - x \ra \leq r^2\}, \quad x \in \dom f,
\ea
$$
belongs to $\dom f$ for any $r \in [0,1)$. Thus, for $\dom f$ containing no straight line, the Hessian is positive definite at any feasible point.

In this paper, we often use the local norms defined by the Hessians of SCF. For $x \in \dom f$, $s \in \dom f_*$, $u \in \E$, and $g \in \E^*$, we adopt the following notation:
$$
\ba{rclrcl}
\| u \|_{\nabla^2 f(x)} & = & \la \nabla^2 f(x) u, u \ra^{1/2}, &  \quad \| g \|_{\nabla^2 f(x)} & = & \la  g, [\nabla^2 f(x)]^{-1} g \ra^{1/2}, \\
\\
\| g \|_{\nabla^2 f_*(s)} & = & \la g, \nabla^2 f_*(s) g \ra^{1/2}, &  \quad \| u \|_{\nabla^2 f_*(s)} & = & \la [\nabla^2 f_*(s)]^{-1} u, u \ra^{1/2}.
\ea
$$
If no ambiguity arise, we use the corresponding shortcuts $\| u \|_x$, $\| g \|_x$, $\| g \|_s$, and $\| u \|_s$. Hence, the sense of notation $\| a \|_b$ depends on the sets and the spaces containing $a$ and $b$. 

For SCF, we often use the following relations:
\beq\label{eq-Hess}
\ba{rcl}
(1-r)^2 \nabla^2 f(x) & \preceq & \nabla^2 f(u) \; \preceq \; {1 \over (1-r)^2} \nabla^2 f(x), \quad x \in \dom f,
\ea
\eeq
which are valid for all $y \in W^f_r(x)$ and $r \in [0,1)$. The following inequalities are also important:
\beq\label{eq-SCF}
\ba{rcl}
\omega(r) & \leq & f(y) - f(x) - \la \nabla f(x), y - x \ra \; \leq \; \omega_*(r),
\ea
\eeq
where $x \in \dom f$, $y \in W^f_r(x)$, $r \in [0,1)$, and
$$
\ba{rcl}
\omega(\tau) & = & \tau - \ln(1+\tau), \quad \omega_*(\tau) \; = \; - \tau - \ln(1-\tau), \quad \tau \in [0,1).
\ea
$$
Sometimes the following variant of inequalities (\ref{eq-SCF}) is useful:
\beq\label{eq-SCF1}
\ba{rcl}
\omega(\delta) \; \leq \; f(x) - f(y) - \la \nabla f(y), x - y \ra & \leq &  \omega_*(\delta),
\ea
\eeq
with $\delta = \| \nabla f(y) - \nabla f(x) \|_x$ (for (\ref{eq-SCF1}$_2$), we need $\delta < 1$, see Theorem 5.1.12 in \cite{Lect}).

The following fact is often absent in the textbooks.
\BL\label{lm-DGrad}
Let function $f(\cdot)$ be self-concordant. Then, for any $x \in \dom f$ and $y \in \E$ with $r = \| y - x \|_x < 1$, we have
\beq\label{eq-DGrad}
\ba{rcl}
\| \nabla f(y) - \nabla f(x) - \nabla^2 f(x)(y-x) \|_x & \leq & {r^2 \over 1 - r}.
\ea
\eeq
\EL
\proof
Note that $\Delta_g \Def \nabla f(y) - \nabla f(x) - \nabla^2 f(x)(y-x) = G(y-x)$ with
$$
\ba{rcl}
G & = & \int\limits_0^1 [\nabla^2 f(x+\tau(y-x)) - \nabla^2 f(x)]d \tau \\
\\
& \stackrel{(\ref{eq-Hess})}{\preceq} & \nabla^2 f(x) \int\limits_0^1 \Big[{1 \over (1- \tau r)^2} - 1 \Big] d \tau \; = \; {r \over 1 - r}  \nabla^2 f(x),
\ea
$$
and $G \stackrel{(\ref{eq-Hess})}{\succeq} \nabla^2 f(x) \int\limits_0^1 [(1- \tau r)^2 - 1] d \tau \; \succeq \; -{r \over 1 - r}  \nabla^2 f(x)$. Therefore,
$$
\ba{rcl}
\| \Delta_g \|^2_x & = & \la G(y-x), [\nabla^2 f(x) ]^{-1} G(y-x) \ra \; \leq {r^2 \over (1-r)^2} \la \nabla^2 f(x)(y-x), y - x \ra. \QF
\ea
$$

Let us recall now the properties of {\em self-concordant barriers } (SCB) for convex cones.
\BD\label{def-SCB}
A self-concordant function $F(\cdot)$ is called a {\em self-concordant barrier}, if there exists a constant 
\beq\label{eq-DefNu}
\ba{rcl}
\nu & \geq & 1
\ea
\eeq
such that for all $x \in \dom F$ and $h \in \E$ we have
\beq\label{eq-SCB}
\ba{rcl}
\la \nabla F(x), h \ra^2 & \leq & \nu \la \nabla^2 F(x) h, h \ra.
\ea
\eeq
The constant $\nu$ is called the {\em parameter} of the barrier.
If $\nabla^2 F(x) \succ 0$, then the condition (\ref{eq-SCB}) is equivalent to the following:
\beq\label{def-SCB1}
\ba{rcl}
\lambda_F^2(x) \; \Def \; \| \nabla F(x) \|_x^2 & \leq & \nu, \quad x \in \dom F.
\ea
\eeq
\ED

This value is responsible for complexity of the set $\dom F$ for corresponding IPMs. A {\em lower bounds} for $\nu$ can be obtained by the following statement (see Theorem 5.4.1 in \cite{Lect}).
\BT\label{th-LowB}
Let $p_i$, $i=1, \dots , k$, be recession directions of set $Q$ and $\bar x \in \inter Q$. Define coefficients $\{\beta_i\}_{i=1}^k$ such that
$$
\ba{rcl}
\bar x - \beta_i p_i & \not\in & \inter Q, \quad i = 1, \dots, k.
\ea
$$
And let $\bar x - \sum\limits_{i=1}^k \alpha_i p_i \in Q$ for some positive $\alpha_i$, $i = 1 \dots, k$. Then $\nu \geq \sum\limits_{i=1}^k {\alpha_i \over \beta_i}$.
\ET

In the theory of IPMs, the most important feasible sets are convex cones. A cone $K$ is called {\em proper} if it is closed convex and pointed (contains no straight lines). For such cones, SCBs possess a natural property of {\em logarithmic homogeneity}:
\beq\label{def-HomB}
\ba{rcl}
F(\tau x) & \equiv & F(x) - \nu \ln \tau, \quad x \in \inter K, \; \tau > 0.
\ea
\eeq
This identity has several important consequences: for all $x \in \inter K$ and $\tau > 0$, we have
\beq\label{eq-Hom12}
\ba{rclrcl}
\nabla F(\tau x) & = & {1 \over \tau} \nabla F(x), & \nabla^2 F(\tau x) & = & {1 \over \tau^2} \nabla^2 F(x),
\ea
\eeq
\beq\label{eq-GX}
\ba{rcl}
\la \nabla F(x), x \ra & = & - \nu, \quad \la \nabla^2 F(x) x, x \ra \; = \; \nu,
\ea
\eeq
\beq\label{eq-HomHX}
\ba{rclrcl}
\nabla^2 F(x) x & = & \nabla F(x), &  \quad \| \nabla F(x) \|^2_x & = & \nu.
\ea
\eeq
Thus, the parameter of a logarithmically homogeneous SCB for a proper cone (which we call {\em regular barrier}) is equal to the degree of logarithmic homogeneity.

It is important that the dual barrier $F_*(\cdot)$ for the regular barrier $F(\cdot)$ is a regular barrier for the dual cone
$$
\ba{rcl}
K^* & = & \Big\{ s \in \E^*: \; \la s, x \ra \geq 0, \; \forall x \in K \Big\},
\ea
$$
with the same value of barrier parameter. Note that for  logarithmically homogeneous barriers, we can use a stronger version of Fenchel inequality:
\beq\label{eq-HFenchel}
\ba{rcl}
F(x) + F_*(s) & \geq &  \nu \ln {\la s, x \ra \over \nu} - \nu, \quad x \in \inter K, \; s \in \inter K^*,
\ea
\eeq
where the equality is achieved only for $s = - \mu \nabla F(x)$ with some $\mu > 0$.

\BE\label{ex-BarSDP}
Let $\S^n = \{X \in \R^{n \times n}:\; X = X^T  \}$. Denote $\S^n_+ \Def \{ X \in \S^n: \; X \succeq 0\}$. It is a proper cone, which admits the following self-concordant logarithmically homogeneous barrier:
$$
\ba{rcl}
F(X) & = & - \ln \det X, \quad \nu = n.
\ea
$$
Note that $(\S^n_+)^* = \S^n_+$ and $F_*(S) = - \ln \det S - n$. The derivatives of these barriers are very simple:
\beq\label{def-DetDir}
\ba{rcl}
\nabla F(X) & = & - X^{-1}, \quad \nabla^2 F(X) H \; = \; X^{-1} H X^{-1}, \\
\\
\; [\nabla^2F(X)]^{-1}H & = & X H X, \quad X \succ 0, \; H \in \S^n. \QR
\ea
\eeq
\EE

\section{Primal and dual central paths}\label{sc-PDPaths}
\SetEQ

Consider the following pair of primal-dual conic problems:
\beq\label{prob-PD}
\ba{rcl}
f^* \; = \; \min\limits_{{\cal F}_p = \left\{\ba{c} x \in K \\ A x = b \ea \right.} \la c, x \ra & =  & \max\limits_{{\cal F}_d = \left\{ \ba{c} y \in \R^m, s \in K^* \\ s + A^* y = c \ea \right.} \la b, y \ra, \ea
\eeq
where
$K \subset \E$ is a proper cone, $K^* \subset \E^*$ is its dual cone,
$b \in \R^m \setminus \{0\}$, and $A$ is a linear operator from $\E$ to $\R^m$. Thus, $A^*: \R^m \to \E^*$.
We assume that the extended operator $(A^*,c): \R^{m+1} \to \E^*$ has full column rank. We assume also existence of a strictly feasible primal-dual point
$$
\ba{rcl}
(\tilde x, \tilde y, \tilde s) & \in & \inter K \times \R^m \times \inter K^*: \quad A \tilde x = b, \quad s + A^* \tilde y = c.
\ea
$$
This condition is sufficient for existence and boundedness of the primal-dual optimal set and existence of the central path (e.g. \cite{LongStep}).
\BE\label{ex-SDP}
Consider the standard problem of SDO:
\beq\label{prob-SDOP}
\ba{c}
\min\limits_{X \in K} \Big\{ \la C, X \ra: \; \la A_i, X \ra = b^{(i)}, \; i = 1, \dots, m \Big\},
\ea
\eeq
where $C$ and $A_i$, $i = 1, \dots, m$, are symmetric $n \times n$-matrices, and $K = \S^n_+$.
In this formulation, the linear operator $A(\cdot): \S^n \to \R^m$ is defined as follows:
$$
\ba{rcl}
A^{(i)}(X) & = & \la A_i, X \ra, \quad i = 1, \dots, m.
\ea
$$
Therefore, the conjugate operator $A^*(\cdot): \R^m \to \S^n$ acts as $A^*(y) = \sum\limits_{i=1}^m y^{(i)} A_i$, $y \in \R^m$. Since $K = K^* = \S^n_+$, the problem dual to (\ref{prob-SDOP}) is as follows:
\beq\label{prob-SDOD}
\ba{rcl}
\max\limits_{y \in \R^m} \Big\{ \la b, y \ra: \; C \succeq A^*(y) =  \sum\limits_{i=1}^m y^{(i)} A_i \Big\}. \QF
\ea
\eeq
\EE

Let cone $K$ admit a computable $\nu$-normal barrier $F(\cdot)$. Denote by $F_*(\cdot)$ its dual barrier, which we assume to be computable too. However, it may happen that the computational cost of the primal barrier is higher than that one of the dual barrier.

The standard approach for treating the pair of problems (\ref{prob-PD}) consists in unifying them in a single {\em primal-dual problem} with zero optimal value:
\beq\label{prob-UPD}
\ba{rcl}
\min\limits_{z=(x,y,s)} \Big\{ \la c, x \ra - \la b, y \ra : \; x \in {\cal F}_p, \; (y,s) \in {\cal F}_d \Big\} \; = \; 0.
\ea
\eeq
This problem can be solved by a path-following scheme, based on the primal-dual barrier function $\Phi(z) = F(x) + F_*(s)$, $z \in {\cal F} \Def {\cal F}_p \times {\cal F}_d$.  The majority of the efficient interior-point methods follow somehow the primal-dual central path defined as follows:
\beq\label{def-PDCP}
\ba{rcl}
z(t) & \Def & \arg\min\limits_{z \in {\cal F}} \Big\{ F(x) + F_*(s): \; \la c, x \ra - \la b,y \ra = {\nu \over t} \Big\}, \quad t>0.
\ea
\eeq
This process is usually divided on two stages.
\begin{itemize}
\item
{\em Corrector Stage}. At this stage, we fix some value of penalty parameter $t > 0$ and try to find a point $z \approx z(t)$ by solving approximately the minimization problem~(\ref{def-PDCP}). The closeness of point $z$ to the points at the central path can be estimated by different proximity measures. For our purposes, the most convenient is the following {\em Functional Proximity Measure} \cite{LongStep}:
\beq\label{def-FProx}
\ba{rcl}
\Omega(z) & = & \nu \ln \la s, x \ra + F(x) + F_*(s) + \nu - \nu \ln \nu\\
\\
& = & \nu \ln [\la c, x \ra - \la b, y \ra] + F(x) + F_*(s) + \nu - \nu \ln \nu \; \geq \; 0, \quad z \in {\cal F}.
\ea
\eeq
\item
{\em Predictor Stage}. If $z$ is close enough to the central path, we compute direction $d \approx z'(t)$ and move along it, keeping $z + \alpha d \approx z(t + \alpha)$, $\alpha > 0$. It is convenient to estimate this closeness by checking the values $\Omega(z+\alpha d)$.
\end{itemize}

The bigger the value of $\alpha$ is, the higher is efficiency of the scheme. It mainly depends on the quality of approximation $d$ to the exact tangent direction to the central path. For an important family of {\em symmetric cones}, covering the problems of Linear, Quadratic, and Semidefinite Optimization, a good direction $d$ can be defined by a {\em Scaling Point}. 

For these cones, two components of a strictly feasible primal-dual pair $(x,s) \in \inter (K \times K^*)$ can be linked by a uniquely defined scaling point $w \in \inter K$, ensuring the following relation:
\beq\label{eq-ScaleP}
\ba{rcl}
s & = & \nabla^2 F(w) x.
\ea
\eeq
It appears that the Hessian $\nabla^2 F(w)$ can be used for defining a very good approximation $d=(\Delta x, \Delta y, \Delta s)$ to derivative of the central path, which is called {\em Affine-Scaling Directon} (ASD). This direction is a solution of the following system of linear equations:
\beq\label{def-ASD}
\ba{rcl}
\Delta s + \nabla^2F(w) \Delta x & = & p, \quad A \Delta x \; = \; 0, \quad \Delta s + A^* \Delta y \; = \; 0,
\ea
\eeq
where the right-hand side $p \in \E^*$ depends on our particular goals. For example, for defining the {\em centering direction}, we choose $p = -s - \mu \nabla F(x)$ with $\mu = {1 \over \nu} \la s, x \ra$, and for ASD, we take $p = -s$.

Thus, for symmetric cones, the sequence of our actions is as follows.
\begin{enumerate}
\item
Find $z=(x,y,s) \approx z(t)$ ensuring the relation (\ref{eq-ScaleP}).
\item
Compute the scaling point $w \in \inter K$.
\item
Compute the Affine-Scaling Direction by (\ref{def-ASD}).
\end{enumerate}

However, for one of the most important applications of this theory, the Semidefinite Optimization (SDO), the second action is quite expensive. It needs computation of square roots of two symmetric positive-definite matrices, while all other computational tasks are handled by the standard Cholesky factorization.

In this paper, we suggest to use an alternative sequence of actions.
\begin{enumerate}
\item
Compute the scaling point $w \in \inter K$ using an auxiliary minimization process either in the primal or in the dual spaces.
\item
Using the Hessian $\nabla^2F(w)$, form the strictly feasible point $z=(x,y,s) \approx z(t)$, which satisfies the scaling equation (\ref{eq-ScaleP}). 
\item
Compute the Affine-Scaling Direction by (\ref{def-ASD}) with $p = -s$.
\end{enumerate}

As we will see later, with this strategy, all necessary computations  be accomplished by the standard Cholesky factorization. A primal variant of this approach was already discussed in \cite{TW}. In this paper, we justify its dual variant, which is more appropriate for many practical applications.

For our purposes, it is important to define separately the central paths for the primal and for the dual problem:
\beq\label{def-xys}
\ba{rcl}
x_t & = & \arg\min\limits_{A x = b} [ \; f_t(x) \Def t \la c, x \ra + F(x)\; ],\\
\\
(y_t,s_t) & = & \arg\max\limits_{s+A^*y = c} [\; \vf_t(y,s) \Def t \la b, y \ra - F_*(s) \;].
\ea
\eeq
It appears, that this definition provides us with the same trajectories as (\ref{def-PDCP}). Denote 
$$
\ba{rcl}
f^*_t & = & f_t(x_t), \quad \vf^*_t \; = \; \vf_t(y_t,s_t), 
\ea
$$
and define the dual barrier function
\beq\label{def-psi}
\ba{rcl}
\zeta(y) & \Def & F_*(c - A^* y), \quad y \in {\cal F}_y \Def \{ y \in \R^m: \; c - A^* y \in K^* \},\\
\\
\psi_t(y) & = & \zeta(y) - t \la b, y \ra, \quad t > 0.
\ea
\eeq
Note that $\zeta(\cdot)$ is a self-concordant barrier for the set ${\cal F}_y$ with parameter $\nu_{\zeta} \leq \nu$. Our approach is very efficient for the problems with small values of $\nu_{\zeta}$. In Section \ref{sc-Easy}, we discuss examples of SDO-problems with $\nu_{\zeta} = O(m)$. 

\BL\label{lm-PDCP}
For all $t > 0$, we have $z(t) \equiv (x_t, y_t, s_t)$.
Moreover, 
\beq\label{eq-SX}
\ba{rcl}
s_t & = & - {1 \over t} \nabla F(x_t), \quad \la c, x_t \ra - \la b, y_t \ra \; = \; {\nu \over t},
\ea
\eeq
\beq\label{eq-SXF}
\ba{rcl}
F(x_t) + F_*(s_t) & = & -\nu \ln t - \nu,
\ea
\eeq
and
\beq\label{def-FTStar}
\ba{rcl}
f^*_t & = & \vf^*_t + \nu \ln t,
\ea
\eeq
\EL
\proof
The first-order optimality condition for point $x_t$ looks as follows:
\beq\label{eq-LUOpt}
\ba{rcl}
t c + \nabla F(x_t) & = & t A^* u_t,
\ea
\eeq
where $u_t$ is a dual multiplier for the linear equality constraints, which for convenience is multiplied by a constant factor $t$. Then, $c - A^* u_t = - {1 \over t} \nabla F(x_t) \in \inter K^*$, and we have
$$
\ba{rcl}
\nabla \psi_t(u_t) & \refEQ{eq-LUOpt} & - A \nabla F_*(-{1 \over t} \nabla F(x_t)) - t b \refEQ{eq-Hom12} - t A \nabla F_*(- \nabla F(x_t)) - t b \\
\\
& \refEQ{eq-FDual} & t (Ax_t-b) = 0.
\ea
$$
Therefore, $u_t = y_t$, and we get the first equality in (\ref{eq-SX}). It implies
$$
\ba{rcl}
\la c, x_t \ra - \la b, y_t \ra & = & \la c, x_t \ra - \la A x_t, y_t \ra \; \stackrel{(\ref{eq-SX}_1)}{=} \; \la - {1 \over t} \nabla F(x_t), x_t \ra \; \refEQ{eq-HomHX} \; {\nu \over t},
\ea
$$
and the second relation in (\ref{eq-SX}) is proved. Now, we see that
$$
\ba{rcl}
F(x_t) + F_*(s_t) & \stackrel{(\ref{eq-HFenchel}),(\ref{eq-SX}_1)}{=} &
\nu \ln {\la s_t, x_t \ra \over \nu} - \nu \; \refEQ{eq-SX} \; - \nu \ln t - \nu, 
\ea
$$
and this is (\ref{eq-SXF}).

Thus, $(x_t,y_t,s_t)$ is feasible for problem (\ref{def-PDCP}), and for all $z \in {\cal F}$ with $\la c, x \ra - \la b, y \ra = {\nu \over t}$, we have 
$$
\ba{rcl}
F(x) + F_*(s) & \refGE{eq-HFenchel} &
\nu \ln {\la s, x \ra \over \nu} - \nu \; = \; \nu \ln {\la c, x \ra - \la b, y \ra \over \nu} - \nu \\
\\
& = & - \nu \ln t - \nu \; \refEQ{eq-SXF} \; 
\; F(x_t) + F_*(s_t).
\ea
$$
In view of uniqueness of the solution of problem (\ref{def-PDCP}), this means that $z(t) \equiv z_t$.
Finally,
$$
\ba{rcl}
f^*_t - \vf^*_t & = & F(x_t) + F_*(s_t) + t [\la c, x_t \ra - \la b, y_t \ra ] \; \refEQ{eq-SX} \; F(x_t) + F_*(- {1 \over t} \nabla F(x_t)) + \nu \\
\\
& \refEQ{eq-Hom12} & F(x_t) + F_*(- \nabla F(x_t)) + \nu + \nu \ln t \; \stackrel{(\ref{eq-HFenchel}),(\ref{eq-GX})}{=} \; \nu \ln t. \QR
\ea
$$

For $z = (x,y,s) \in {\cal F}$ and $t > 0$, let us define now a {\em partial} proximity measure:
\beq\label{def-OmegaT}
\ba{rcl}
\Omega_t(z) & = & [f_{t}(x) - f^*_{t}] + [\vf^*_{t} - \vf_t(y,s)]\\
\\
& = & t[\la c, x \ra - \la b, y \ra] + F(x) + F_*(s) - f^*_{t} + \vf^*_{t}.
\ea
\eeq
Thus, 
\beq\label{def-OT}
\ba{rcl}
\Omega_t(z) & = & t[\la c, x \ra - \la b, y \ra] + F(x) + F_*(s) - \nu \ln t.
\ea
\eeq

\BL\label{lm-OT}
For any $z \in \rint {\cal F}$, we have $\Omega(z) = \min\limits_{t>0} \Omega_t(z)$. The minimum is attained at 
\beq\label{def-TZ}
\ba{rcl}
t(z) & \Def & {\nu \over \la c, x \ra - \la b, y \ra} = {\nu \over \la s, x \ra}. \QF
\ea
\eeq
\EL

\BC\label{cor-OT}
For $t = t(z)$, we have $\Omega_t(z) = \Omega(z)$. \qed
\EC

For proper interpretation of the subsequent results, it is important to know that the primal and the dual problems in (\ref{prob-PD}) have {\em complementary complexities}. In order to show this, let us look at the barrier function for the primal feasible set
$$
\ba{rcl}
g(x) & = & F(x), \quad x \in \rint {\cal F}_p.
\ea
$$
This is a self-concordant barrier with parameter $\nu_g \leq \nu$.
In order to highlight its relation with the parameter $\nu_{\zeta}$ of the dual barrier $\zeta(\cdot)$, let us prove the following identity.
\BL\label{lm-LNorms}
For any $t>0$, we have
\beq\label{eq-LNorms}
\ba{rcl}
\lambda_g^2(x_t) + \lambda_{\zeta}^2(y_t) & = & \nu.
\ea
\eeq
\EL
\proof
Indeed, for any $x \in \rint {\cal F}_p$, we have
$$
\ba{rcl}
- \half \lambda^2_g(x) & = & \min\limits_{h \in \E} \Big\{ \la \nabla F(x), h \ra + \half \| h \|_x^2 \; : Ah = 0 \Big\}.
\ea
$$
Denote by $h^*$ the optimal solution of this problem and by $u^* \in \R^m$ the vector of optimal dual multipliers. Then
$$
\ba{rcl}
\nabla F(x) + \nabla^2 F(x) h^* & = & A^T u^*, \quad A h^* \; = \; 0.
\ea
$$
Hence, $h^* = [\nabla^2 F(x)]^{-1} \Big[ A^T u^* - \nabla F(x) \Big ] \refEQI{eq-HomHX}{1} [\nabla^2 F(x)]^{-1} A^T u^* + x$. Therefore,
$$
\ba{rcl}
u^* & = & - \Big[ A[\nabla^2 F(x)]^{-1} A^T\Big]^{-1} A x \; \refEQI{prob-PD}{1} \; - \Big[ A[\nabla^2 F(x)]^{-1} A^T \Big]^{-1}b,
\ea
$$
and we conclude that
$$
\ba{rcl} 
\la \nabla F(x), h^* \ra & = & \la \nabla F(x), [\nabla^2 F(x)]^{-1} A^T u^* + x \ra \; \stackrel{(\ref{eq-GX}_1),(\ref{eq-HomHX}_1)}{=} \; - \la x, A^T u^* \ra - \nu,\\ \\
\la \nabla^2F(x) h^*, h^* \ra & = & \la \nabla^2F(x) ([\nabla^2 F(x)]^{-1} A^T u^* + x), [\nabla^2 F(x)]^{-1} A^T u^* + x \ra\\
\\
& \refEQI{eq-GX}{2} & \la \Big[ A[\nabla^2 F(x)]^{-1} A^T \Big]^{-1} b, b \ra + 2 \la x, A^T u^* \ra + \nu.
\ea
$$
Hence, $- \half \lambda^2_g(x) = - \half \nu + \half \la \Big[ A[\nabla^2 F(x)]^{-1} A^T \Big]^{-1} b, b \ra$. For $x_t \refEQI{eq-SX}{1}  - {1 \over t} \nabla F_*(s_t)$, we have 
$$
\ba{rcl}
\nabla^2 F(x_t) & \refEQI{eq-Hom12}{2} & t^2 \nabla^2 F(-\nabla F_*(s_t)) \; \refEQI{eq-FDual}{4} \; t^2 [\nabla^2 F_*(s_t)]^{-1}.
\ea
$$
Therefore,
$$
\ba{rcl}
\nu & = & \lambda^2_g(x_t) + \la \Big[ A[\nabla^2 F(x_t)]^{-1} A^T \Big]^{-1} b, b \ra\; = \; \lambda^2_g(x_t) + t^2\la \Big[ A \nabla^2 F_*(s_t) A^T \Big]^{-1} b, b \ra\\
\\
& \refEQI{def-psi}{1} & \lambda^2_g(x_t) + t^2\la [ \nabla^2 \zeta(y_t)]^{-1} b, b \ra \; \refEQI{def-xys}{2} \; \lambda^2_g(x_t) + \la [ \nabla^2 \zeta(y_t)]^{-1} \nabla \zeta(y_t), \nabla \zeta(y_t) \ra \\
\\
& \refEQ{def-SCB1} & \lambda^2_g(x_t) + \lambda^2_{\zeta}(y_t). \QR
\ea
$$

To the best of our knowledge, equality (\ref{eq-LNorms}) is new. It has several important consequences. First of all, it provides us with strong arguments {\em against} the short-step primal-dual methods. In these methods, the steps are confined by the corresponding Dikin ellipsoids, and the penalty parameters for the primal and the dual central paths are the same. From the theory of self-concordant functions (see Section 5.3.1 in \cite{Self1}), we know that the increase of the current value of penalty parameter $t$ must be in inverse proportion to the local norm of the gradient of corresponding self-concordant barrier. Hence, by (\ref{eq-LNorms}), either the primal or the dual norm is big. This fact cancels any hope for simultaneous large steps in the primal and the dual spaces.

Further, from inequality (\ref{eq-LNorms}), we get
\beq\label{eq-Ass}
\ba{rcl}
\nu & \refLE{def-SCB1} & \nu_g + \nu_{\zeta}.
\ea
\eeq
This relation is another motivation for a proper treating of the {\em Primal-Dual Asymmetry}. Indeed, if the dual feasible set ${\cal F}_y$ is simple (and the barrier parameter $\nu_{\zeta}$ is small), then the primal feasible set ${\cal F}_p$ {\em cannot} be simple since
$$
\ba{rcl}
\nu_g & \refGE{eq-Ass} & \nu - \nu_{\zeta}.
\ea
$$
In this case, it is better to apply strategies related to the dual central path. We will discuss several examples in Section \ref{sc-Easy}. Below, we consider a simple SDO-problem with this type of asymmetry.
\BE\label{ex-SDO1}
Let us choose some $a, b \in \R^n$ with $\la a, b \ra > 0$. Consider the following SDO-problem:
\beq\label{prob-SDO1}
\ba{rcl}
f^* & = & \min\limits_{X \succeq 0} \Big\{ \la I_n, X \ra: \; X a = b \Big\}.
\ea
\eeq
In this problem, we have ${n(n+1) \over 2}$ variables and $n$ linear equality constraints. It is easy to see that its optimal solution is $X^* = {bb^T \over \la a, b \ra}$ with $f^* = {\| b \|^2 \over \la a, b \ra}$.

For rewriting the problem (\ref{prob-SDO1}) in the standard form (\ref{prob-SDOP}), denote $A_i = \half (a e_i^T + e_i a^T)$, $i = 1, \dots, n$. Then the equality constraints  in (\ref{prob-SDO1}) becomes as follows:
$$
\ba{rcl}
\la A_i , X \ra & = & \half \la a e_i^T + e_i a^T, X \ra \; = \; \la e_i, X a \ra \; = \; b^{(i)}, \quad i=1,\dots, n.
\ea
$$
Thus, in accordance to (\ref{prob-SDOD}), the problem dual to (\ref{prob-SDO1}) is as follows:
\beq\label{prob-SDO1D}
\ba{rcl}
\max\limits_{y \in \R^n} \Big\{ \la b, y \ra: \; I_n \succeq \half(y a^T + a y^T)\Big\},
\ea
\eeq
with the barrier function $\zeta(y) = - \ln \det\left(I_n - \half(y a^T + a y^T)\right)$. Let us compute this determinant in the closed firm.

It is easy to see that the matrix $S(y) \Def I_n - \half(y a^T + a y^T)$ has only two eigenvalues different from one:
$$
\ba{rcl}
\lambda_{1,2} & = & 1 - \half \la a,y \ra \pm \half \| a \| \cdot \| y\|.
\ea
$$
Hence, $\det S(y) = (1 - \half \la a,y \ra)^2 - {1 \over 4} \| a \|^2 \| y \|^2$.

Let us compute the barrier parameter $\nu_{\zeta}$. Without loss of generality, we can assume that $a = \alpha e_1$ with $\alpha = \| a \|>0$. Hence,
$$
\ba{rcl}
\det S(y) & = & \Big(1 - \half \alpha y^{(1)}\Big)^2 - {1 \over 4} \alpha^2 \| y \|^2 \; = \; 1 - \alpha y^{(1)} - {1 \over 4} \alpha^2 \sum\limits_{i=2}^n (y^{(i)})^2.
\ea
$$
This is a concave quadratic function and therefore $\nu_{\zeta} = 1$ (see Example 5.3.1(4) in \cite{Lect}). Hence, the restriction $g(\cdot)$ of the barrier $F(X) = - \ln \det X$ onto the feasible set of the primal problem (\ref{prob-SDO1}) has the parameter
\beq\label{eq-ParSDO1}
\ba{rcl}
\nu_{g} & \refGE{eq-Ass} & n - 1. 
\ea
\eeq

Inequality (\ref{eq-ParSDO1}) is valid only for the barrier parameter of our particular barrier function. Let us show that the primal feasible set of problem (\ref{prob-SDO1}) indeed requires a self-concordant barrier with large value of $\nu_g$.

Without loss of generality, we can think that  $a = \alpha e_1$ and $b^{(1)} > 0$. Then the matrix $\bar X = I_n - e_1 e_1^T + X^*$ belongs to the relative interior of the feasible set ${\cal F}_p$ of problem (\ref{prob-SDO1}). Note that the directions $p_i = e_i e_i^T$, $i=2, \dots, n$ are recession directions for ${\cal F}_p$. On the other hand,
$$
\ba{rcl}
\bar X - p_i & \not\in & \rint {\cal F}_p, \quad i=2, \dots , n, \quad \bar X - \sum\limits_{i=2}^n p_i \; = \; X^* \; \in {\cal F}_p.
\ea
$$
Hence, by Theorem \ref{th-LowB}, parameter of any barrier for ${\cal F}_p$ satisfies inequality (\ref{eq-ParSDO1}).
\qed
\EE

This example demonstrates the maximal asymmetry in the complexities of the primal and the dual barrier functions. In such a case, we should work with problem with the minimal value of the barrier parameter. However, in the real-life applications, it is usually impossible to get a priori good bounds for their values. Fortunately, in our situation we discover a small value of barrier parameter for the {\em standard} barrier function. We do not need to change it. As we will see later, an appropriate predictor-corrector strategy indeed can have global complexity bounds automatically dependent on $\min\{ \nu_g, \nu_{\zeta} \}$.

To conclude this section, let us give a simple sufficient condition for relation $\nu_g \geq \nu_{\zeta}$.
\BL\label{lm-FY}
Let set ${\cal F}_y$ be bounded. Then $\nu_g = \nu \geq \nu_{\zeta}$.
\EL
\proof
Indeed, if the set ${\cal F}_y$ is bounded, then there exists its analytic center $y^*_{\zeta}$ with $\nabla \zeta(y^*_{\zeta}) = 0$. Since $\lim\limits_{t \to 0} y_t = y^*_{\zeta}$, we get
$$
\ba{rcl}
\nu \; \geq \; \nu_g & \refGE{def-SCB1} & \lim\limits_{t \to 0} \lambda^2_g(x_t) \; \refEQ{eq-LNorms} \; \lim\limits_{t \to 0}\, [\nu - \lambda^2_{\zeta}(y_t)] \; = \; \nu.
\ea
$$
It remains to note that parameter $\nu_{\zeta}$ does not exceed $\nu$ too.
\qed

Of course, the symmetric statement is also valid: boundedness of ${\cal F}_p$ implies $\nu_{\zeta} = \nu$. Sometimes, this distinction can be seen from the analytical structure of the problem. However, it is not always the case: in Example \ref{ex-SDO1} both ${\cal F}_p$ and ${\cal F}_y$ are unbounded.

\section{Dual Gambit Rule}\label{sc-BDS}
\SetEQ

Let us discuss now a simple way for computing the scaling point. For that, we consider a {\em pure dual} centering process, which finds an approximate solution of the dual problem
\beq\label{prob-DCP}
\ba{c}
\min\limits_{y \in {\cal F}_y} \Big\{ \psi_t(y) = \zeta(y) - t \la b, y \ra \Big\}.
\ea
\eeq
This can be done either by the rudimentary Damped Newton Method:
\beq\label{met-DN}
\ba{rcl}
\lambda_k & = & \la \nabla \psi_t(y_k), [\nabla^2 \zeta(y_k)]^{-1} \nabla \psi_t(y_k) \ra^{1/2}, \\
\\
y_{k+1} & = & y_k - {1 \over 1 + \lambda_k} [\nabla^2 \zeta_d(y_k)]^{-1} \nabla \psi_t(y_k), \quad k \geq 0,
\ea
\eeq
or by an auxiliary path-following scheme (e.g. Section 5.3.5 in \cite{Lect}). Recall that
\beq\label{eq-fdHess}
\ba{rcl}
\nabla \zeta(y) & = & - A \nabla F_*(c- A^*y), \quad \nabla^2 \zeta(y) \; = \; A \nabla^2 F_*(c-A^*y) A^*.
\ea
\eeq

In any case, we can assume that after a finite number of iterations, we get a point $\bar y$, satisfying the following condition 
\beq\label{eq-DTerm}
\ba{rcl}
\bar \lambda \Def \la \nabla \psi_t(\bar y), [\nabla^2 \zeta(\bar y)]^{-1} \nabla \psi_t(\bar y) \ra^{1/2} & \leq & \beta \; < \; 1.
\ea
\eeq
In our computations, we use the Newton direction 
\beq\label{def-DBar}
\ba{rcl}
\bar d & \Def & [ \nabla^2 \zeta(\bar y)]^{-1} \nabla \psi_t(\bar y),
\ea
\eeq
which is also required for checking the condition (\ref{eq-DTerm}). Note that
\beq\label{eq-BDNorm}
\ba{rcl}\
\| \bar d \|^2_{\bar y} & \refEQ{eq-DTerm} & \la \nabla^2 \zeta(\bar y) \bar d, \bar d \ra \; = \; \bar \lambda^2.
\ea
\eeq

Denote $\bar s = c-A^*\bar y \in \inter K^*$ and $\bar x  = - {1 \over t} \nabla F_*(\bar s) \in \inter K$. Then
\beq\label{eq-BarSX}
\ba{rcl}
\bar s & \stackrel{(\ref{eq-FDual}),(\ref{eq-Hom12})}{=} & - {1 \over t} \nabla F(\bar x). 
\ea
\eeq

Consider now a special rule for generating a feasible point in a close neighbourhood of the primal-dual central path. 
\beq\label{eq-Out}
\ba{|c|}
\hline \\
\mbox{\sc Dual Gambit Rule}\\
\\
\hline \\
\ba{rclrcl}
\hat y & = & \bar y + \bar d, & \hat s & = & c - A^* \hat y \; \equiv \; \bar s - A^* \bar d \;\\
\\
\bar w_* & = & \sqrt{t} \, \bar s, & \hat x & = & \nabla^2 F_*(\bar w_*)\hat s.
\ea\\ \\
\hline
\ea
\eeq
\BE\label{eq-OutSDO}
Let us show how it works for the dual SDO-problem (\ref{prob-SDOD}). In this case, 
$$
\ba{rcl}
\zeta(y) & = & - \ln \det \left(C - \sum\limits_{i=1}^m y^{(i)} A_i \right).
\ea
$$
Let the point $\bar y$ be an approximate solution of the problem (\ref{prob-DCP}). After computing $\bar d$ by definition (\ref{def-DBar}), we form $\hat y = \bar y + d$ and compute $\hat S = C - \sum\limits_{i=1}^m \hat y^{(i)} A_i$. Then 
\beq\label{eq-XOut}
\ba{rcl}
\hat X & \refEQI{def-DetDir}{2} & {1 \over t} \bar S^{-1} \hat S \bar S^{-1}. 
\ea
\eeq
As we will see later, the popular rule (\ref{eq-ScaleSDP}) for computing firstly the scaling point is much more expensive.
\qed
\EE

Note that the point $\hat y$ is computed in (\ref{eq-Out}) by the {\em backward} Newton step from the point $\bar y$, which necessarily {\em increases} the value of the dual function $\psi_t(\cdot)$:
$$
\ba{rcl}
\psi_t(\hat y) & \geq & \psi_t(\bar y) + \la \nabla \psi_t(\bar y), \hat y - \bar y \ra \; \refEQ{eq-BDNorm} \; \psi_t(\bar y) + \bar \lambda^2.
\ea
$$
However, it appears that point $\hat z = (\hat x, \hat y, \hat s)$ is feasible and it is close to the central path.
\BL\label{lm-Out}
We have 
\beq\label{eq-HXFeas}
\ba{rcl}
A \hat x & = & b,
\ea
\eeq
and
\beq\label{eq-HDist}
\ba{rcl}
\| \hat x - \bar x \|_{\bar x} & = & \| \hat s - \bar s \|_{\bar s} \; = \; \bar \lambda \leq \beta < 1. 
\ea
\eeq
Hence, $\hat x \in \inter K$, $\hat s \in \inter K^*$, and
$\hat z = (\hat x,\hat y,\hat s) \in \rint {\cal F}$.
\EL
\proof
Note that $\| \hat s - \bar s \|_{\bar s} = \| A^* \bar d \|_{\bar s} \refEQ{eq-fdHess} \bar \lambda < 1$. Hence, $\hat s \in \inter K^*$. 
Further,
$$
\ba{rcl}
0 & = & \nabla^2 \psi_t(\bar y) \bar d - \nabla \psi_t(\bar y) \; \refEQ{eq-fdHess} \; A \Big[ \nabla^2 F_*(\bar s) A^* \bar d + \nabla F_*(\bar s)\Big]+ t b \\
\\
& \stackrel{(\ref{eq-HomHX}_1)}{=} &  A \nabla^2 F_*(\bar s) [A^* d  - \bar s] + t b \; \stackrel{(\ref{eq-Hom12}_2)}{=} \; 
t \Big[ A \nabla^2 F_*(\bar w_*) [A^* d  - \bar s] + b\Big] \; \refEQ{eq-Out} \; t[b - A \hat x].
\ea
$$

Finally, in view of (\ref{eq-Hom12}$_2$), we have $t \nabla^2 F_*(\bar \omega_*) = \nabla^2 F_*(\bar s)$. Therefore,
\beq\label{eq-DT}
\ba{rcl}
t \hat x - t \bar x & \refEQ{eq-Out} & \nabla^2F_*(\bar s) (\bar s - A^*\bar d) + \nabla F_*(\bar s) \; \stackrel{(\ref{eq-HomHX}_1)}{=} \; - \nabla^2 F_*(\bar s) A^* \bar d. 
\ea
\eeq
Hence,
$$
\ba{rcl}
\| \hat x - \bar x \|^2_{\bar x} & = & {1 \over t^2} \| \nabla^2 F_*(\bar s) A^* \bar d \|^2_{\bar x} \; \stackrel{(\ref{eq-Hom12}_2)}{=} \; \la \nabla^2 F(t\bar x)\nabla^2 F_*(\bar s) A^* d, \nabla^2 F_*(\bar s) A^* d \ra\\
\\
& \stackrel{(\ref{eq-FDual}_2)}{=} &\la \nabla^2 F_*(\bar s) A^* d, A^* d \ra \; \refEQ{eq-BDNorm} \; \bar \lambda^2 \; \leq \beta^2 < 1.
\ea
$$
Thus, $\hat x \in \inter K$, and we conclude that $\hat z \in \rint {\cal F}$.
\qed

Let us present now a geometric interpretation of the point $\hat x$. 
\BL\label{lm-XPlus}
We have
\beq\label{eq-XPFirst}
\ba{rcl}
t c + \nabla F(\bar x) + \nabla^2 F(\bar x)(\hat x - \bar x) & = & t A^*(\bar y - d).
\ea
\eeq
Hence, 
\beq\label{eq-XProj}
\ba{rcl}
\hat x & = & \arg\min\limits_{A x = b} \Big\{ \la \nabla f_t(\bar x), x - \bar x \ra + \half \| x - \bar x \|^2_{\bar x} \Big\}.
\ea
\eeq 
\EL
\proof
Indeed,
$$
\ba{rcl}
t c + \nabla F(\bar x) + \nabla^2 F(\bar x)(\hat x - \bar x) & \refEQ{eq-DT} & t c - t \bar s + t^2 [\nabla^2 F_*(\bar s)]^{-1} \Big( - {1 \over t} \nabla^2 F_*(\bar s)A^* \bar d \Big)\\
\\
& = & t c - t \bar s - t A^*\bar d \; = \; t A^* (\bar y - \bar d).
\ea
$$
This relation, in view of equality $A \hat x = b$, justifies that $\hat x$ is an optimal solution of problem~(\ref{eq-XProj}).
\qed

It appears that the point $\hat x$ is very close to the primal central path.
\BL\label{lm-LF}
$\lambda_{f_t}(\hat x) \leq \left({\beta \over 1 - \beta} \right)^2$.
\EL
\proof
Indeed, 
$$
\ba{rcl}
\lambda^2_{f_t}(\hat x) & = & \max\limits_{A h = 0} \Big\{ 2 \la t c + \nabla F(\hat x), h \ra - \| h \|^2_{\hat x} \Big\}\\
\\
& \refEQ{eq-XPFirst} & \max\limits_{A h = 0} \Big\{ 2 \la \nabla F(\hat x) - \nabla F(\bar x) - \nabla^2F(\bar x)(\hat x - \bar x), h \ra - \| h \|^2_{\hat x} \Big\}\\
\\
& \leq & \| \nabla F(\hat x) - \nabla F(\bar x) - \nabla^2F(\bar x)(\hat x - \bar x) \|^2_{\hat x}\\
\\
& \refLE{eq-HDist} & { 1 \over (1 - \beta)^2}\| \nabla F(\hat x) - \nabla F(\bar x) - \nabla^2F(\bar x)(\hat x - \bar x) \|^2_{\bar x} \; \refLE{eq-DGrad} \; {\beta^4 \over (1-\beta)^4}. \QF
\ea
$$

Note that 
$$
\ba{rcl}
\psi_t(\hat y) \refLEI{eq-SCF}{2}\psi_t(\bar y) + \bar \lambda^2 + \omega_*(\bar \lambda) & \refLE{eq-SCF1} & \psi^*_t + \bar \lambda^2 + 2 \omega_*(\bar \lambda). 
\ea
$$
Hence, the point $\hat y$ remains in a close neighbourhood of the dual central path. This is also confirmed by the following statement.
\BL\label{lm-LPsi}
We have $\lambda_{\psi_t}(\hat y) \leq {2 \beta - \beta^2\over (1 - \beta)^2}$.
\EL
\proof 
Indeed, 
$\lambda_{\psi_t}(\hat y) = \| \nabla \psi_t(\hat y) \|^*_{\hat y} \leq {1 \over 1-\beta} \| \nabla \psi_t(\hat y) \|^*_{\bar y}$. It remains to note that
$$
\ba{rcl}
\| \nabla \psi_t(\hat y) \|^*_{\bar y} & \leq & \| \nabla \psi_t(\hat y) - \nabla \psi_t(\bar y) - \nabla^2 \zeta(\bar y)(\hat y - \bar y) \|^*_{\bar y} + \|\nabla \psi_t(\bar y) + \nabla^2 \zeta(\bar y) \bar d \|^*_{\bar y}\\
\\
& \refLE{eq-DGrad} & {\beta^2 \over 1 - \beta} + 2 \beta \; = \; {2 \beta - \beta^2 \over 1 - \beta}. \QR
\ea
$$

We can see that a small step back in the quality of the dual point $\bar y$ is compensated by a much higher quality of the point $\hat x$ and by  existence of the scaling point $\bar w_*$. This is the reason for calling the rule (\ref{eq-Out}) of forming a strictly feasible primal-dual solution the {\em Dual Gambit Rule} (DG-rule).

Thus, we conclude that the point $\hat z$ is close to the primal-dual central path:\footnote{In the last inequality, we use the following simple fact on univariate convex functions: if $\xi(\cdot)$ is convex and $\xi(0)=0$, then for all $a,b \geq 0$, we have $\xi(a)+\xi(b) \leq \xi(a+b)$.}
\beq\label{eq-Close}
\ba{rcl}
\Omega(\hat z) & \leq & \Omega_t(\hat z) \; = \; [f_t(\hat x) - f^*_t] + [\vf_t^* - \vf_t(\hat y)]\\
\\
& \stackrel{{\bf L.} \ref{lm-LF}, {\bf L.}\ref{lm-LPsi}}{\leq} & \omega_*\left( {\beta^2 \over (1-\beta)^2} \right) + \omega_*\left({2 \beta - \beta^2 \over (1- \beta)^2} \right) \; \leq \; \omega_*\left({2 \beta \over (1- \beta)^2} \right).
\ea
\eeq

Finally, we need an estimate for the closeness of the reference value $t(\hat z)$ to $t$.
\BL\label{lm-DifT}
We have $\left(1 - {\beta \over \sqrt{\nu}} \right)^2 \leq {t \over t(\hat z)} \leq \left(1 + {\beta \over \sqrt{\nu}} \right)^2$.
\EL
\proof
Indeed, $ {t\nu \over t(\hat z)} = t \la \hat s, \hat x \ra \refEQ{eq-Out} t \la \hat s, \nabla^2 F_*(\bar w_*) \hat s \ra \refEQ{eq-Hom12} \| \hat s \|^2_{\bar s}$. It remains to note that 
$$
\ba{rcl}
\| \hat s \|_{\bar s} & \refLE{eq-HDist} & \| \bar s \|_{\bar s} + \beta \refEQ{eq-HomHX} \sqrt{\nu} + \beta.
\ea
$$
Similarly, we can show that $\| \hat s \|_{\bar s} \geq \sqrt{\nu} - \beta$.
\qed

\section{Affine-Scaling Direction}\label{sc-ASD}
\SetEQ

For analyzing the properties of Affine-Scaling Direction, we need several auxiliary statements.
Let us introduce the primal scaling point 
\beq\label{def-SPrim}
\ba{rcl}
\bar w & \Def & - \nabla F_*(\bar w_*) \refEQ{eq-Out} - \nabla F_*(\sqrt{t} \bar s) \refEQI{eq-Hom12}{1} - {1 \over \sqrt{t} } \nabla F_*(\bar s)\\
\\
& \refEQ{eq-BarSX} & - {1 \over \sqrt{t} } \nabla F_*(- {1 \over t} \nabla F(\bar x)) \; \refEQI{eq-Hom12}{1} \; - \sqrt{t} \, \nabla F_*(-\nabla F(\bar x)) \; \refEQ{eq-FDual} \;\sqrt{t} \bar x \in \inter K.
\ea
\eeq
Then,
\beq\label{eq-HPrim}
\ba{rcl}
\nabla^2 F(\bar w) & \refEQI{eq-Hom12}{2} &  {1 \over t} \nabla^2 F(\bar x) \; \refEQ{eq-BarSX} \; {1 \over t} \nabla^2 F(- {1 \over t} \nabla F_*(\bar s)) \; \refEQI{eq-Hom12}{2} \; t \nabla^2 F(-  \nabla F_*(\bar s))\\
\\
& \refEQ{eq-FDual} & \; t [\nabla^2 F_*(\bar s)]^{-1} \; \refEQI{eq-Hom12}{2} \; [\nabla^2 F_*(\sqrt{t} \bar s)]^{-1} \; \refEQ{eq-Out} \; [\nabla^2 F_*(\bar w_*)]^{-1}.
\ea
\eeq

Using the point $\bar w$, the scaling equation in (\ref{eq-Out}) can be rewritten as follows:
\beq\label{eq-PScale}
\ba{rcl}
\hat s & = & \nabla^2 F(\bar w) \hat x.
\ea
\eeq

For self-scaled cones \cite{Self1}, equation (\ref{eq-PScale}) implies also 
\beq\label{eq-GScale}
\ba{rcl}
\nabla F(\hat x) & = & \nabla^2 F(\bar w) \nabla F_*(\hat s). 
\ea
\eeq 
However, for the general cones, we can prove only a relaxed version of this equality. We need to justify first the following fact.
\BL\label{lm-Relax}
Let point $\hat z$ be formed by the rule (\ref{eq-Out}). Then, for any $\tau \in \R$, we have
\beq\label{eq-VarTau}
\ba{rcl}
\nabla F(\bar x) + \tau \nabla^2 F(\bar x) (\hat x - \bar x) & = & \nabla^2 F(\bar w) \Big[ \nabla F_*(\bar s) + \tau \nabla^2 F_*(\bar s) (\hat s - \bar s)\Big].
\ea
\eeq
At the same time,
\beq\label{eq-ScaleG}
\ba{rcl}
\| \nabla F_*(\hat s) - [\nabla^2 F(\bar w)]^{-1} \nabla F(\hat x)\|_{\bar w} & \leq & {2\beta^2 \over 1 - \beta} \sqrt{t}.
\ea
\eeq
\EL
\proof
Since $\nabla^2 F(\bar w) \nabla F_*(\bar s) \refEQ{eq-BarSX} - t \nabla^2 F(\bar w) \bar x \refEQI{eq-Hom12}{2} - \nabla^2 F(\bar x) \bar x \refEQI{eq-HomHX}{1}  \nabla F(\bar x)$, for proving (\ref{eq-VarTau}), it is enough to apply the following transformations:
$$
\ba{rcl}
\nabla^2 F(\bar w) \nabla^2 F_*(\bar s) \hat s & \refEQI{eq-Hom12}{2} & t \nabla^2 F(\bar w) \nabla^2 F_*(\bar w_*) \hat s \; \refEQ{eq-HPrim} \; t \hat s \\
\\
& \refEQ{eq-PScale} & t \nabla^2 F(\bar w) \hat x \; \refEQ{def-SPrim} \; \nabla^2 F(\bar x) \hat x.
\ea
$$ 

Further, denoting $h_p = \hat x - \bar x$ and $h_d = \hat s - \bar s$, and applying (\ref{eq-VarTau}) with $\tau = 1$, we can see that
$$
\ba{c}
\| \nabla F_*(\hat s) - [\nabla^2 F(\bar w)]^{-1} \nabla F(\hat x)\|_{\bar w}\\
\\
= \| \nabla F_*(\hat s)  - \nabla F_*(\bar s) - \nabla^2 F_*(\bar s) h_d - [\nabla^2 F(\bar w)]^{-1} [ \nabla F(\hat x) - \nabla F(\bar x) - \nabla^2 F(\bar x) h_p] \|_{\bar w}
\\
\\
\leq \| \delta_a \|_{\bar w}+ \| [\nabla^2F(\bar w)]^{-1} \delta_b \|_{\bar w} = \| \delta_a \|_{\bar w}+ \| \delta_b \|_{\bar w}, \\
\\
\delta_a \Def \nabla F_*(\hat s) - \nabla F_*(\bar s) - \nabla^2 F_*(\bar s) h_d, \quad 
\delta_b \Def\nabla F(\hat x) - \nabla F(\bar x) - \nabla^2 F(\bar x) h_p.
\ea
$$

Note that $\| h_d \|_{\bar s} \refEQ{eq-HDist} \bar \lambda \leq \beta$. Then
$$
\ba{c}
\| \delta_a \|^2_{\bar w} \; \refEQ{def-SPrim} \; \la \nabla^2 F(\sqrt{t} \bar x) \delta_a, \delta_a \ra \; \refEQI{eq-Hom12}{2} \; {1 \over t} \la \nabla^2 F(\bar x) \delta_a, \delta_a \ra \; \refEQ{eq-BarSX} \;
{1 \over t} \la \nabla^2 F(-{1 \over t} \nabla F_*(\bar s)) \delta_a, \delta_a \ra \\
\\
\refEQI{eq-Hom12}{2} \;  t \la \nabla^2 F(- \nabla F_*(\bar s)) \delta_a, \delta_a \ra\; \refEQ{eq-FDual} \;
t \la [\nabla^2 F_*(\bar s)]^{-1} \delta_a, \delta_a \ra \; \refLE{eq-DGrad} \; t {\beta^4 \over (1- \beta)^2}.
\ea
$$ 
Similarly, $\| h_p \|_{\bar x} \refEQ{eq-HDist} \bar \lambda \leq \beta $. Hence
$$
\ba{rcl}
\| \delta_b \|^2_{\bar w} & \refEQ{def-SPrim} & t \la \delta_b,  [\nabla^2F(\bar x)]^{-1} \delta_b \ra \; \refLE{eq-DGrad} \; t {\beta^4 \over (1-\beta)^2}. \QF
\ea
$$

In the standard theory of interior-point methods for symmetric cones, the linear system~(\ref{def-ASD}) is solved with different right-hand sides, which provide us either with tangent or with centering step. In our approach, the centering procedure is reduced to solving the dual problem~(\ref{prob-DCP}). Hence, we need to define only the tangent direction $\Delta z = (\Delta x, \Delta y, \Delta s)$, which can be computed by the following system of linear equations:
\beq\label{def-TanD}
\ba{rcl}
\Delta s + \nabla^2 F(\bar w) \Delta x & = & - \hat s, \quad A \Delta x \; = \; 0, \quad \Delta s + A^* \Delta y \; = \; 0,
\ea
\eeq
where $\bar w$ is defined by (\ref{def-SPrim}), and $\bar x$, $\bar s$, and $\hat s$ are defined by the DG-rule (\ref{eq-Out}). This system can be also written in the primal form:
\beq\label{def-TanD1}
\ba{rcl}
\Delta x + \nabla^2 F_*(\bar w_*) \Delta s & = & - \hat x, \quad A \Delta x \; = \; 0, \quad \Delta s + A^* \Delta y \; = \; 0.
\ea
\eeq

We need the following properties of this direction.
\BL\label{lm-AffS}
We have the following relations:
\beq\label{eq-AOrt}
\ba{rcl}
\la \Delta s, \Delta x \ra & = & 0,
\ea
\eeq
\beq\label{eq-AChange}
\ba{rcl}
\la \hat s, \Delta x \ra + \la \Delta s, \hat x \ra & = & - \la \hat s, \hat x \ra,
\ea
\eeq
\beq\label{eq-AGap}
\ba{rcl}
\la c, \hat x + \Delta x \ra - \la b, \hat y + \Delta y \ra & = & 0,
\ea
\eeq
\beq\label{eq-ANorm}
\ba{rcl}
\| \Delta x \|^2_{\bar w} + \| \Delta s \|^2_{\bar w} & = & \la \hat s, \hat x \ra.
\ea
\eeq
\EL
\proof
Indeed, $\la \Delta s, \Delta x \ra \refEQI{def-TanD}{3} - \la A^* \Delta y, \Delta x \ra = \la A \Delta x, \Delta y \ra \refEQI{def-TanD}{2} 0$. Further,
$$
\ba{rcl}
\la \hat s, \Delta x \ra + \la \Delta s, \hat x \ra & \refEQ{def-TanD} & \la \hat s, \Delta x \ra - \la \hat s + \nabla^2 F(\bar w) \Delta x, \hat x \ra \; \refEQ{eq-BarSX} \; - \la \hat s, \hat x \ra.
\ea
$$
This is equation (\ref{eq-AChange}).
At the same time, we have
$$
\ba{c}
\la c, \hat x + \Delta x \ra - \la b, \hat y + \Delta y \ra \; \refEQ{eq-HXFeas} \; 
\la \hat s + A^* \hat y, \hat x + \Delta x \ra - \la A \hat x, \hat y + \Delta y \ra\\
\\
\refEQI{def-TanD}{2} \; \la \hat s, \hat x + \Delta x \ra - \la A \hat x, \Delta y \ra \; \refEQ{eq-AChange} \; - \la \Delta s, \hat x \ra - \la A^*\Delta y, \hat x \ra \; \refEQI{def-TanD}{3} \; 0.
\ea
$$
Finally, in view of (\ref{eq-AOrt}), we have
$$
\ba{rcl}
\| \Delta x \|^2_{\bar w} & \refEQI{def-TanD}{1} & - \la \hat s, \Delta x \ra, \quad \| \Delta s\|^2_{\bar w} \; \refEQI{def-TanD1}{1} \; - \la \Delta s, \hat x\ra.
\ea
$$
It remains to use (\ref{eq-AChange}).
\qed

Let us bound the growth of primal-dual proximity measure along the affine-scaling direction $\Delta z$. 
\BT\label{th-Predict}
For the points $z(\alpha) = \hat z + \alpha \Delta z$, we have the following upper bound for the growth of the proximity measure:
\beq\label{eq-PGrow}
\ba{rcl}
\Omega(z(\alpha)) & \leq & \Omega(\hat z)+ \alpha c_{\beta} \beta^2  + \omega_*(\alpha c_{\beta}) - \nu \omega_*(\alpha)\\
\\
& \leq & \Omega(\hat z) - \ln (1 - \alpha c_{\beta}),
\ea
\eeq
where $c_{\beta} = {\sqrt{\nu}+ \beta \over 1-\beta}$.
\ET
\proof
Note that $\ln {\la \hat s + \alpha \Delta s, \hat x + \alpha \Delta x \ra \over \la s, x \ra} \refEQ{eq-AOrt} \ln \left(1 + {\alpha \over  \la \hat s, \hat x \ra} [\la \hat s , \Delta x \ra + \la \Delta s, \hat x \ra ] \right) \refEQ{eq-AChange} \ln(1-\alpha)$. Therefore,
$$
\ba{rcl}
\Omega(z(\alpha) - \Omega(\hat z) & = & F(\hat x + \alpha \Delta x) + F_*(\hat s + \alpha \Delta s) + \nu \ln {\la \hat s + \alpha \Delta s, \hat x + \alpha \Delta x \ra \over \la \hat s, \hat x \ra}  - F(\hat x) - F_*(\hat s)\\
\\
& \refLE{eq-SCF} & \alpha \Big[ \la \nabla F(\hat x), \Delta x \ra + \la \Delta s, \nabla F_*(\hat s) \ra \Big]_A + \omega_*(\alpha \rho) + \nu \ln (1-\alpha),
\ea
$$
where $\rho^2 = \la \nabla^2 F(\hat x) \Delta x, \Delta x \ra + \la \Delta s, \nabla^2 F_*(\hat s) \Delta s \ra$. In view of (\ref{eq-HDist}), we have
$$
\ba{rcl}
\rho^2 & \refLE{eq-Hess} & {1 \over (1-\beta)^2} \Big[ \la \nabla^2 F(\bar x) \Delta x, \Delta x \ra + \la \Delta s, \nabla^2 F_*(\bar s) \Delta s \ra \Big] \\
\\
& \refEQI{eq-Hom12}{2} & {t \over (1-\beta)^2} \Big[ \la \nabla^2 F(\bar w) \Delta x, \Delta x \ra + \la \Delta s, \nabla^2 F_*(\bar w_*) \Delta s \ra \Big] \\
\\
& \refEQ{eq-ANorm} & {t \la \hat s, \hat x \ra \over (1-\beta)^2} \; \refEQ{def-TZ} \; {t \over t(\hat z)} \cdot {\nu \over (1-\beta)^2} \; \stackrel{\rm \bf L.\ref{lm-DifT}}{\leq} \; \left( {\sqrt{\nu}+ \beta \over 1-\beta} \right)^2 \Def \; c_{\beta}^2.
\ea
$$
On the other hand,
$$
\ba{rcl}
\la \Delta s + \nabla^2 F(\bar w) \Delta x, \nabla F_*(\hat s) \ra & \refEQ{def-TanD} &  \la - \hat s, , \nabla F_*(\hat s) \ra \; \refEQ{eq-GX} \; \nu,\\
\\
\la \nabla F(\hat x), \Delta x + [\nabla^2 F(\bar w) ]^{-1} \Delta s \ra & \refEQ{def-TanD1} & - \la \nabla F(\hat x), \hat x \ra \; \refEQ{eq-GX} \; \nu.
\ea
$$
Therefore,
$$
\ba{rcl}
a & \Def & \la \nabla^2 F(\bar w)\Delta x - \Delta s, \nabla F_*(\hat s) - [\nabla^2 F(\bar w) ]^{-1} \nabla F(\hat x) \ra\\
\\
& = &   \la \nabla^2 F(\bar w)\Delta x, \nabla F_*(\hat s) \ra +  \la \Delta s, [\nabla^2 F(\bar w) ]^{-1} \nabla F(\hat x) \ra - [ \cdot ]_A\\
\\
& = & \la \nabla^2 F(\bar w)\Delta x + \Delta s, \nabla F_*(\hat s) \ra +  \la \nabla F(\hat x) , [\nabla^2 F(\bar w) ]^{-1} \Delta s + \Delta x\ra - 2[ \cdot ]_A\\
\\
& \stackrel{(\ref{def-TanD}),(\ref{def-TanD1})}{=} & \la - \hat s, \nabla F_*(\hat s) \ra +  \la \nabla F(\hat x) , - \hat x\ra - 2[ \cdot ]_A
\; \refEQ{eq-GX} \; 2 \nu - 2 [ \cdot ]_A.
\ea
$$
Thus, $[\cdot ]_A = \nu - \half a$. It remains to find an upper bound for $|a|$. 
$$
\ba{rcl}
|a| & \leq & \| \nabla^2 F(\bar w)\Delta x - \Delta s \|_{\bar w} \cdot \| \nabla F_*(\hat s) - [\nabla^2 F(\bar w) ]^{-1} \nabla F(\hat x) \|_{\bar w}\\
\\
& \refEQ{eq-AOrt} & \Big[ \| \Delta x \|^2_{\bar w} + \| \Delta s \|_{\bar w}^2 \Big]^{1/2} \cdot \| \nabla F_*(\hat s) - [\nabla^2 F(\bar w) ]^{-1} \nabla F(\hat x) \|_{\bar w}\\
\\
& \refLE{eq-ScaleG} & \Big[ \| \Delta x \|^2_{\bar w} + \| \Delta s \|_{\bar w}^2 \Big]^{1/2} \cdot {2 \beta^2 \over 1 - \beta} \sqrt{t}\\
\\
& \refEQ{eq-ANorm} & {2 \beta^2 \over 1 - \beta} \sqrt{t \la \hat s, \hat x \ra} \; = \; {2 \beta^2 \sqrt{\nu} \over 1 - \beta} \sqrt{t \over t(\hat z)} \; \stackrel{\rm \bf L. \ref{lm-DifT}}{\leq} { 2\beta^2(\sqrt{\nu} + \beta) \over 1 - \beta} \; = \; 2 \beta^2 c_{\beta}.
\ea
$$
Thus, we come to the first inequality in (\ref{eq-PGrow}). The second one is valid since $\beta \in [0,1)$.
\qed

\BR\label{rm-PGrow}
For self-scaled cones, in the above proof we have $a \refEQ{eq-GScale} 0$. In this case, the growth of the proximity measure can be bounded as follows:
\beq\label{eq-PGrowSS}
\ba{rcl}
\Omega(z(\alpha)) & \leq & \Omega(\hat z) + \omega_*(\alpha c_{\beta}) - \nu \omega_*(\alpha).
\ea
\eeq
\ER

\section{General dual primal-dual predictor-corrector \\ method}\label{sc-DBMet}
\SetEQ

Let us assume that the set ${\cal F}_y= \{ y \in \R^n: \; c - A^*y \in K^* \}$ is bounded. In accordance to Lemma \ref{lm-FY}, this means that $\nu_{\zeta} \leq \nu_g = \nu$. Hence, it is natural to associate both the centering and the predicting strategies with the dual problem.

In view of our assumption, there exists an {\em analytic center} $y^*_{\zeta}$ of the set ${\cal F}_y$, characterized by the first-order optimality condition
\beq\label{def-ACenter}
\ba{rcl}
\nabla \zeta(y^*_{\zeta}) & = & 0,
\ea
\eeq
(see Section 5.3.3 in \cite{Lect}).
With respect to this point, the {\em asphericity} of the dual feasible set ${\cal F}_y$ is bounded by the barrier parameter as follows:
\beq\label{eq-ASF}
\ba{rcl}
W^{\zeta}_1(y^*_{\zeta}) & \subseteq & {\cal F}_y \; \subseteq \; W^{\zeta}_{\rho_y}(y^*_{\zeta}),  \quad \rho_y \Def \nu_{\zeta} + 2 \sqrt{\nu_{\zeta}}.
\ea
\eeq

Consider the following switching strategy.
\beq\label{met-DBS}
\ba{|l|}
\hline \\
\hspace{2ex} \mbox{\sc General Predictor-Corrector Method with Dual Gambit Step}\\
\\
\hline \\
\mbox{{\bf Initialization.} Choose $\beta \in (0,{1 \over 4}]$ and define $A = \delta + \omega_*\left({2 \beta \over (1-\beta)^2}\right)$ with $\delta >0$. }\\
\\
\mbox{Compute $y_0 \in {\cal F}_y$ with $\| \nabla \zeta(y_0) \|_{y_0} \leq \half \beta$ and set $t_0 = {\beta \over 2 \| b \|_{y_0}}$.}\\
\\
\mbox{\bf $k$th iteration ($k \geq 0$)}\\ \\
\ba{rl}
\mbox{a)} & \mbox{For $y_k \in \inter {\cal F}_y$ and $t_k > 0$, compute $B_k=[ \nabla^2 \zeta(y_k)]^{-1}$, $g_k = \nabla \psi_{t_k}(y_k)$,}\\
\\
& \mbox{$d_k = B_k g_k$, and $\lambda_k = \la g_k, d_k \ra^{1/2}$. Define $s_k = c - A^* y_k \in \inter K^*$.}\\
\\
\mbox{b)} & \mbox{{\bf If} $\lambda_k > \beta$, {\bf Then} $y_{k+1} = y_k - {d_k \over 1 + \lambda_k}$, $t_{k+1} = t_k$. (Damped Newton Step)}\\
\\
\mbox{c)} & \mbox{{\bf Else} (Predictor Step)}\\ \\
& \mbox{1. Set $w^*_k = t_k^{1/2} s_k$, $\hat y_k = y_k + d_k$, $\hat s^k = c - A^* \hat y_k$, and $\hat x_k = \nabla^2 F_*(w^*_k) \hat s_k$.}\\
\\
& \mbox{2. Find Affine-Scaling Direction $\Delta z_k = (\Delta x_k, \Delta y_k, \Delta s_k)$ by the system}\\
\\
& \hspace{5ex} \mbox{$\Delta x_k + \nabla^2 F_*(w^*_k) \Delta s_k = - \hat x_k, \quad A \Delta x_k = 0, \quad \Delta s_k + A^* \Delta y_k = 0$.}\\
\\
& \mbox{3. Compute the step size $\alpha_k > 0$ from the equation $\Omega(\hat z_k + \alpha_k \Delta z_k) = A$.}\\
\\
& \mbox{4. Define $z_{k+1} = \hat z_k + \alpha_k \Delta z_k$ and $t_{k+1} = t(z_{k+1}) \refEQQ{def-TZ}{eq-AChange} {\nu \over (1-\alpha_k)\la \hat s_k, \hat x_k \ra}$.}
\ea\\
\\
\hline
\ea
\eeq

\BR\label{rm-ASD}
It is important that step {\rm c2)} in (\ref{met-DBS}) does not require additional matrix inversions. Indeed, at this moment, matrix $B_k$ is already computed. At the same time,
$$
\ba{rcl}
0 & = & A \Delta x_k \; \refEQ{def-TanD1} \; A\Big[ - \hat x_k + {1 \over t_k} \nabla^2 F_*(s_k) A^* \Delta y_k \Big] \;  \refEQ{eq-HXFeas} \; {1 \over t_k} \nabla^2 \zeta(y_k) \Delta y_k - b.
\ea
$$
Hence, using the direction $h_k = B_k b$, we can define 
\beq\label{eq-FormASD}
\ba{rcl}
(\Delta y_k,\Delta s_k) & = & t_k (h_k, - A^* h_k), \quad \Delta  x_k \; = \; \nabla^2 F_*(s_k) A^* h_k - \hat x_k. \QF
\ea
\eeq
\ER

For the future use, let us prove the following bound.
\BL\label{lm-GrowDY}
The norm of ASD is bounded as follows:
\beq\label{eq-GrowDY}
\ba{rcl}
\| \Delta y_k \|_{\hat y_k} & \leq & {\beta + \lambda_{\zeta}(y_k) \over 1 - \beta} \; \leq \; \hat c_{\beta} \Def  {\beta + \sqrt{\nu_{\zeta}} \over 1 - \beta}.
\ea
\eeq
\EL
\proof
Note that $\| y_k - \hat y_k \|_{y_k} \refEQI{eq-fdHess}{2} \| A^T d_k \|_{s_k} \refEQI{eq-HDist}{2} \lambda_k \leq \beta$. Therefore,
$$
\ba{rcl}
(1-\beta)\| \Delta y_k \|_{\hat y_k} & \refLE{eq-Hess} & \| \Delta y_k \|_{y_k} \; \refEQ{eq-FormASD} \; \| t_k b \|_{y_k} \; \refLE{def-psi} \; \| \nabla \psi_{t_k}(y_k) \|_{y_k} + \| \nabla \zeta(y_k) \|_{y_k}\\
\\
& \refEQ{met-DBS} & \lambda_k + \lambda_{\zeta}(y_k) \;  \refLE{def-SCB1} \; \beta + \sqrt{\nu_{\zeta}}. \QR
\ea
$$

Let us estimate now the performance of method (\ref{met-DBS}). It alternates two groups of iterations of different type. The first type corresponds to a Correction Process in the dual space, implemented by Step b). The second group, which always consists of one iteration, is the Predictor Step implemented by Step c) in the primal-dual space.

During the Correction Process, the values of penalty parameter $t_k$ remain unchanged. Let us look at the iterations of a correction group, started at iteration $k_0$ and finished at iteration $k_1$. Since $t_k = t_{k_0}$ and $\lambda_k > \beta$ for all $k_0 \leq k \leq k_1-1$, each iteration in this group decreases the value of penalty function $\psi_{t_{k_0}}(\cdot)$ as follows:
\beq\label{eq-FDec}
\ba{rcl}
\psi_{t_{k_0}}(y_{k+1}) & \leq & \psi_{t_{k_0}}(y_{k}) - \omega(\beta).
\ea
\eeq

By the choice of the starting point $y_0$, we have
$$
\ba{rcl}
\lambda_0 & = & \| \nabla \zeta(y_0) - t_0 b \|_{y_0} \; \leq \; \| \nabla \zeta(y_0) \|_{y_0} + t_0 \| b \|_{y_0} \; \leq \; \beta.
\ea
$$
Hence, at iteration $k=0$, we always perform the predictor step.
Therefore, at any iteration $k_0 \geq 1$, the point $z_{k_0}$ is well defined by the Step c$_4$) of the previous predicting iteration. Thus,
$$
\ba{rcl}
\psi_{t_{k_0}}(y_{k_0}) - \psi^*_{t_{k_0}} & \refEQ{def-psi} &
\vf^*_{t_{k_0}} - \vf_{t_{k_0}}(y_{k_0}, c - A^* y_{k_0}) \\
\\
& \refLEI{def-OmegaT}{1} & \Omega_{t_{k_0}}(z_{k_0}) \; \stackrel{\bf L.\ref{lm-OT}}{=} \; \Omega(z_{k_0}) \; \refEQ{met-DBS} \; A.
\ea
$$
Therefore, the length of the Correction Group is bounded:
$k_1 - k_0 \refLE{eq-FDec} {1 \over \omega(\beta)} A$.

For the predictor step at iteration $k_1$, we have
\beq\label{eq-TK1}
\ba{rcl}
t_{k_1+1} & \refEQ{def-TZ} & {t(\hat z_{k_1}) \over 1 - \alpha_k} \; \stackrel{\rm \bf L.\ref{lm-DifT}}{\geq} \; {t_{k_1} \over 1 - \alpha_{k_1}} \left(1 + {\beta \over \sqrt{\nu}} \right)^{-2}  \; \geq \; \exp\left( \alpha_{k_1} - {2 \over \sqrt{\nu}} \beta \right) \cdot t_{k_1}.
\ea
\eeq
Since $\lambda_{k_1} \leq \beta$, we have $\Omega(\hat z_{k_1}) \refLE{eq-Close} \omega_*\left({2 \beta \over (1-\beta)^2}\right)$. Thus,
$$
\ba{rcl}
\delta & \refEQ{met-DBS} & A - \omega_*\left({2 \beta \over (1-\beta)^2}\right) \; \leq \; \Omega(\hat z_{k_1} + \alpha_{k_1} \Delta z_{k_1}) - \Omega(\hat z_{k_1}) \\
\\
& \refLEI{eq-PGrow}{2} & - \ln (1 - \alpha_{k_1} c_{\beta}) \; \leq \; { \alpha_{k_1} c_{\beta} \over 1 -  \alpha_{k_1} c_{\beta}}.
\ea
$$
Hence, 
\beq\label{eq-AGrow}
\ba{rcl}
\alpha_{k_1} - {2 \over \sqrt{\nu}} \beta & \geq & {\delta \over c_{\beta} (1 + \delta)}  - {2\over \sqrt{\nu}}  \beta  \; \refEQ{eq-PGrow} \; {1 \over \sqrt{\nu} + \beta } \left( {\delta (1 - \beta) \over (1 + \delta)}  - 2  \beta \sqrt{1 + {\beta \over \nu}}\right) \\
\\
& \refGE{eq-DefNu} & {1 \over \sqrt{\nu} + \beta } \left( {\delta (1 - \beta) \over (1 + \delta)}  - \beta(2+\beta) \right) \; = \; {1 \over c_{\beta}} \left( {\delta \over 1 + \delta}  - {\beta(2+\beta) \over 1 - \beta} \right).
\ea
\eeq

Thus, we have proved the following theorem.
\BT\label{th-DBS}
Let the paramerts of method (\ref{met-DBS}) satisfy the following conditions:
\beq\label{eq-condDBS}
\ba{rcl}
\beta & \in & (0,{1 \over 4}], \quad \delta \; > \; 0, \quad \vk(\beta,\delta) \; \Def \; {\delta \over 1 + \delta}  - {\beta(2+\beta) \over 1 - \beta} \; > \; 0,
\ea
\eeq
and the starting point $y_0 \in {\cal F}_y$, satisfy the condition
\beq\label{eq-CondY0}
\ba{rcl}
\| \nabla \zeta(y_0) \|_{y_0} & \leq & \half \beta.
\ea
\eeq
Then, each predicting iteration $k \geq 0$ of 
method (\ref{met-DBS}) ensures the following growth of penalty coefficient:
\beq\label{eq-TGrow}
\ba{rcl}
t_0 \; = \; {\beta \over 2 \| b \|_{y_0}}, \quad
t_{k+1} & \geq &  t_k \cdot \exp \Big\{ {k \over \sigma(\beta,\delta)} \Big\}, \quad k \geq 0,
\ea
\eeq
where $\sigma(\beta,\delta) = {\sqrt{\nu} + \beta \over (1-\beta) \vk(\beta,\delta)}\left(\delta + \omega_*\left({2 \beta \over (1-\beta)^2} \right)\right)$.
\ET
\proof
The rate (\ref{eq-TGrow}) follows from inequalities (\ref{eq-TK1}) and (\ref{eq-AGrow}).
\qed

Note that in view of the first inclusion in (\ref{eq-ASF}), we have 
\beq\label{eq-ASF2}
\ba{rcl}
\| b \|_{y^*_{\zeta}} & \leq & f^* - \la b, y^*_{\zeta} \ra.
\ea
\eeq

\BC\label{cor-DBS}
For any $\epsilon > 0$, method (\ref{met-DBS}) generates an $\epsilon$-solution of problem (\ref{prob-PD}) at most in 
\beq\label{eq-CompDBS}
\ba{c}
{\sqrt{\nu} + \beta \over (1-\beta) \vk(\beta,\delta)}\left(\delta + \omega_*\left({2 \beta \over (1-\beta)^2} \right)\right) \cdot \ln\left( {2 \over \epsilon} \cdot [f^* - \la b, y^*_{\zeta} \ra ] \cdot  {(\beta + \sqrt{\nu})^2 \over \beta(1-\beta)} \right)
\ea
\eeq
predicting iterations.
\EC
\proof
Indeed, let at some iteration $k \geq 0$ we have $\lambda_k \leq \beta$. Then the point $\hat z_k$, which is computed at Step c$_1$), is feasible. At the same time,
$$
\ba{rcl}
0 & \leq & \la c, \hat x_k \ra - \la b, \hat y_k \ra \; = \; \la \hat s_k, \hat x_k \ra  \; \refEQQ{eq-Out}{eq-Hom12} \;
{1 \over t_k} \la \hat s_k, \nabla^2F_*(s_k) \hat s_k \ra \\
\\
& \leq & {1 \over t_k} \left( \| \hat s_k - s_k  \|_{s_k} + \| s_k \|_{s_k} \right)^2 \; \refLEQ{eq-HDist}{eq-GX} \; {1 \over t_k} \left( \beta + \sqrt{\nu} \right)^2.
\ea
$$
Thus, for getting an $\epsilon$-solution of problem (\ref{prob-PD}), we need to ensure $t_k \geq {1 \over \epsilon}\left( \beta + \sqrt{\nu} \right)^2$. In view of the lower bound (\ref{eq-TGrow}), it is enough to have 
$$
\ba{rcl}
{\beta \over 2 \| b \|_{y_0}} \exp \Big\{ {k \over \sigma(\beta,\delta)} \Big\} & \geq & {1 \over \epsilon}\left( \beta + \sqrt{\nu} \right)^2.
\ea
$$
At the same time,
$$
\ba{rcl}
\| b \|_{y_0} & \refLE{eq-HDist} & {1 \over 1 - \beta} \| b \|_{y^*_{\zeta}} \; \refLE{eq-ASF2} \; {1 \over 1 - \beta} [f^* - \la b, y^*_{\zeta}\ra].
\ea
$$
It remains to put all inequalities together.
\qed

The complexity bound (\ref{eq-CompDBS}) has one hidden drawback. It is proportional to the square root of the barrier parameter $\nu$. However, the method (\ref{met-DBS}) is essentially dual. Hence, we could expect that the number of iterations of the corresponding interior-point methods is proportional to $\nu_{\zeta}^{1/2}$, which can be much smaller than $\nu^{1/2}$. This is indeed true for an important family of {\em symmetric cones}, which we consider in Section \ref{sc-Sym}. It is important that the functional proximity measure there can be expressed completely in terms of the dual barrier $\zeta(\cdot)$.

\section{Dual primal-dual method for symmetric cones}\label{sc-Sym}
\SetEQ

Let us consider now the problem (\ref{prob-PD}), where the cone $K$ is a regular {\em symmetric cone}. Cone $K$ is symmetric if and only if it admits a normal {\em self-scaled barrier} \cite{Self1}.
\BD\label{def-SSB}
Function $F(\cdot)$ is called a $\nu$-self-scaled barrier for the regular cone $K$ if it is a $\nu$-normal barrier for $K$ and for all $x, w \in \inter K$ we have
\beq\label{eq-IncSSB}
\ba{rcl}
\nabla^2 F(w) x & \in & \inter K^*,\\
\\
F_*(\nabla^2 F(w) x) & = & F(x) - 2 F(w) - \nu. \QF
\ea
\eeq
\ED

It can be proved that for any $w \in \inter K$ we have $K^* = \nabla^2 F(w) K$. For any pair of points $x \in \inter K$, $s \in \inter K^*$ there exists a unique {\em scaling point} $w \in \inter K$ such that
\beq\label{eq-ScaleSWX}
\ba{rcl}
s & = & \nabla^2 F(w)x.
\ea
\eeq
Moreover, the gradients and the Hessians at these points are related as follows:
\beq\label{eq-MapGH}
\ba{rcl}
\nabla F(x) & = & \nabla^2 F(w) \nabla F_*(s), \\
\\
\nabla^2 F(x) & = & \nabla^2 F(w) \nabla^2 F_*(s) \nabla^2 F(w).
\ea
\eeq
\BE\label{eq-SDP}
Let $K \equiv \S^n_+$ be the cone of positive-semidefinite $n \times n$-matrices. Then the following barrier is self-scaled \cite{Self1}:
$$
\ba{rl}
F(X) = & - \ln \det X, \quad F_*(S) \; = \; - \ln \det S - n, \; X \in \inter K, \; S \in \inter K_* = \S^n_+.
\ea
$$
For two positive-definite matrices $S$ and $X$, the scaling point $W$ can be computed by the following expression:
\beq\label{eq-ScaleSDP}
\ba{rcl}
W & = & X^{1/2} [X^{1/2} S X^{1/2}]^{-1/2}X^{1/2}. \QF
\ea
\eeq
\EE

Note that the dual function $F_*(\cdot)$ is a self-scaled barrier for the cone $K^*$. Since the relation (\ref{eq-ScaleSWX}) can be rewritten in the dual form
\beq\label{eq-ScaleXWS}
\ba{rcl}
x & = & \nabla^2 F_*(w_*) s, \quad w_* \Def - \nabla F(w) \in \inter K^*,
\ea
\eeq
we have
\beq\label{eq-RepSSB}
\ba{rcl}
F([\nabla^2 F(w)]^{-1} s) & \refEQ{eq-FDual} & 
F(\nabla^2 F_*(w_*) s) \; = \; F_*(s) - 2 F_*(w_*) - \nu\\
\\
& \refEQI{eq-GX}{1} & F_*(s) + 2 F(w) + \nu.
\ea
\eeq

The above properties are responsible for a very convenient representation of the proximity measure $\Omega(\cdot)$ along ASD.
\BL
Let cone $K$ be self-scaled and the points $x, w \in \inter K$ and $s \in \inter K^*$ are related by (\ref{eq-ScaleSWX}). Then for affine-scaling direction $\Delta z = (\Delta x, \Delta y, \Delta s)$ defined by equations
\beq\label{def-DirSWX}
\ba{rcl}
\Delta s + \nabla^2F(w) \Delta x & = & - s, \quad A \Delta x \; = \; 0, \quad \Delta s + A^* \Delta y \; = \; 0,
\ea
\eeq
we have
\beq\label{eq-OGrow}
\ba{rcl}
\xi_z(\alpha) \; \Def \; \Omega(z + \alpha \Delta z) - \Omega(z) & = & F_*(s + \alpha \Delta s) + F_*\left(s - {\alpha \over 1 - \alpha} \Delta s \right) - 2F_*(s)\\
\\
& = & F(x + \alpha \Delta x) + F\left(x - {\alpha \over 1 - \alpha} \Delta x \right) - 2F(x).
\ea
\eeq
\EL
\proof
Indeed,
$$
\ba{rcl}
\xi_z(\alpha) & \refEQ{def-FProx} & \nu \ln {\la s + \alpha \Delta s, x + \alpha \Delta x \ra \over \la s, x \ra} + F(x+\alpha \Delta x) + F_*(s+\alpha \Delta s) - F(x) - F_*(s)\\
\\
& \refEQQ{eq-AOrt}{eq-AChange} & \nu \ln (1-\alpha) + F(x+\alpha \Delta x) + F_*(s+\alpha \Delta s) - F(x) - F_*(s)\\
\\
& \refEQI{eq-IncSSB}{2} & \nu \ln (1-\alpha) + F(x+\alpha \Delta x) + F_*(s+\alpha \Delta s) - 2 F_*(s) - 2 F(w) - \nu.
\ea
$$
Note that
$$
\ba{rcl}
F(x+\alpha \Delta x) & \refEQ{def-DirSWX} & F(x - \alpha [\nabla^2F(w)]^{-1}(s + \Delta s)) \refEQ{eq-ScaleSWX} \; F((1-\alpha)x - \alpha [\nabla^2F(w)]^{-1} \Delta s)\\
\\
& \refEQ{eq-ScaleSWX} & F( [\nabla^2F(w)]^{-1}((1-\alpha) s - \alpha \Delta s)) \\
\\
& \refEQ{eq-RepSSB} & F_*((1-\alpha)s - \alpha \Delta s) + 2 F(w) + \nu.
\ea
$$
Hence, by logarithmic homogeneity of the barrier $F_*(\cdot)$, we obtain the first representation in (\ref{eq-OGrow}). The second one  can be justified in a similar way.
\qed

\BC\label{cor-RepOG}
In terms of the objects generated by method (\ref{met-DBS}), we have the following representation for the growth of the functional proximity measure:
\beq\label{eq-OGrow1}
\ba{rcl}
\xi_{\hat z_k}(\alpha) & = & \zeta(\hat y_k + \alpha \Delta y_k) + \zeta \left( \hat y_k - {\alpha \over 1 - \alpha} \Delta y_k \right) - 2 \zeta(\hat y_k)\\
\\
& = & F(\hat x_k + \alpha \Delta x_k) + F \left( \hat x_k - {\alpha \over 1 - \alpha} \Delta x_k \right) - 2 F(\hat x_k)
\ea
\eeq
where $\Delta y_k = t_k [\nabla^2 \zeta(y_k)]^{-1} b$ and $\Delta x_k = \nabla^2 F_*(s_k) A^* h_k - \hat x_k$.
\EC
\proof It is enough to use definition (\ref{def-psi}) of the barrier function and the relation (\ref{eq-FormASD}).
\qed

Note that both representations in (\ref{eq-OGrow1}) are valid for the same {\em primal-dual} proximity measure $\Omega(\cdot)$. The first one involves only elements of the dual problem. Hence, let us estimate its growth in terms of the dual barrier function $\zeta(\cdot)$ as follows:
$$
\ba{rcl}
\xi_{\hat z_k}(\alpha) & \refLE{eq-SCF} & \la \nabla \zeta(\hat y_k), \alpha \Delta y_k - {\alpha \over 1 - \alpha} \Delta y_k \ra + \omega_*(\alpha \hat c_{\beta}) + \omega_*\left({\alpha \over 1 - \alpha} \hat c_{\beta} \right)\\
\\
& \refLE{eq-GrowDY} & {\alpha^2 \over 1 - \alpha} \lambda_{\zeta}(\hat y_k) \| \Delta y_k \|_{\hat y_k} + \omega_* \left( {\alpha \hat c_{\beta}(2 - \alpha) \over 1 - \alpha} \right)\\
\\
& \refLEQ{def-SCB1}{eq-GrowDY} & {\alpha^2 \over 1 - \alpha} \sqrt{\nu_{\zeta}} \hat c_{\beta} + \omega_* \left( {2\alpha \hat c_{\beta}\over 1 - \alpha} \right) \; \leq \; {1-\beta  \over 1 - \alpha}(\alpha \hat c_{\beta})^2 + \omega_* \left( {2\alpha \hat c_{\beta}\over 1 - \alpha} \right).
\ea
$$

As in the proof of Theorem \ref{th-DBS}, we conclude that the step size $\alpha_{k_1}$ at the predictor iteration $k_1$ satisfies inequality
\beq\label{eq-DProg}
\ba{rcl}
\delta & \leq & \xi_{\hat z_{k_1}}(\alpha_{k_1})  \; \leq \; {1-\beta  \over 1 - \alpha_{k_1}}(\alpha_{k_1} \hat c_{\beta})^2 + \omega_* \left( {2\alpha_{k_1} \hat c_{\beta}\over 1 - \alpha_{k_1}} \right)
\ea
\eeq
Let us assume that $\alpha_{k_1} \hat c_{\beta} \leq {2 \over 1 - \beta}$. Then
${1-\beta  \over 1 - \alpha_{k_1}}(\alpha_{k_1} \hat c_{\beta})^2 \leq {2\alpha_{k_1} \hat c_{\beta}\over 1 - \alpha_{k_1}}$, and we get
$$
\ba{rcl}
\delta & \refLE{eq-DProg} & - \ln \left(1 - {2\alpha_{k_1} \hat c_{\beta}\over 1 - \alpha_{k_1}} \right) \; \leq \; {{2\alpha_{k_1} \hat c_{\beta}\over 1 - \alpha_{k_1}} \over 1  - {2\alpha_{k_1} \hat c_{\beta}\over 1 - \alpha_{k_1}}}.
\ea
$$
Hence, ${2\alpha_{k_1} \hat c_{\beta}\over 1 - \alpha_{k_1}} \geq {\delta \over 1 + \delta}$. This means that ${\alpha_{k_1} \over 1 - \alpha_{k_1}} \geq {\delta \over 2 \hat c_{\beta}(1 + \delta)}$, and therefore
$\alpha_{k_1}  \geq {\delta  \over 2 \hat c_{\beta} (1+\delta) + \delta}$.
Combining this observation with our temporary assumption, we conclude that 
$$
\ba{rcl}
\alpha_{k_1}  & \geq & \min\Big\{ {2 \over (1-\beta) \hat c_{\beta}}, {\delta  \over 2 \hat c_{\beta} (1+\delta) + \delta} \Big\} \; = \; {\delta  \over 2 \hat c_{\beta} (1+\delta) + \delta} \; \geq \; {\delta  \over (2 \hat c_{\beta}+1) (1+\delta)}.
\ea
$$
Hence, as in the proof of Theorem \ref{th-DBS}, we get the lower bound for the size of predictor step:
$$
\ba{rcl}
\alpha_{k_1} - {2 \over \sqrt{\nu}} \beta & \geq & {1 \over 2 \hat c_{\beta}+1} \left( {\delta  \over 1+\delta} - {2 \beta (2 \hat c_{\beta}+1) \over \sqrt{\nu_{\zeta}}} \right) \; \refGE{eq-DefNu} \; {1 \over 2 \hat c_{\beta}+1} \left( {\delta  \over 1+\delta} - {2 \beta (3+\beta) \over 1 - \beta} \right).
\ea
$$
Now we can present the main theorem on the performance of method (\ref{met-DBS}) for symmetric cones. Its justification is the same as that one for Theorem \ref{th-DBS}.

\BT\label{th-DBS1}
Let cone $K$ in problem (\ref{prob-PD}) be symmetric and
parameters of method (\ref{met-DBS}) satisfy the following conditions:
\beq\label{eq-condDBS1}
\ba{rcl}
\beta & \in & (0,1), \quad \delta \; > \; 0, \quad \hat \vk(\beta,\delta) \; \Def \; {\delta \over 1 + \delta}  - {2\beta(3+\beta) \over 1 - \beta} \; > \; 0.
\ea
\eeq
If the starting point $y_0 \in {\cal F}_y$ satisfies  condition (\ref{eq-CondY0}), then any predictor iteration $k \geq 0$ of
method (\ref{met-DBS}) ensures the following growth of penalty coefficient:
\beq\label{eq-TGrow1}
\ba{rcl}
t_ 0 \; \geq \; {\beta \over 2 \| b \|_{y_0}} , \quad 
t_{k+1} & \geq & t_k \cdot \exp \Big\{ {k \over \hat \sigma(\beta,\delta)} \Big\},
\ea
\eeq
where $\hat \sigma(\beta,\delta) = {2 \hat c_{\beta}+ 1 \over \hat \vk(\beta,\delta)}\left(\delta + \omega_*\left({2 \beta \over (1-\beta)^2} \right)\right)$ and $\hat c_{\beta} = {\beta + \sqrt{\nu_{\zeta}} \over 1 - \beta}$.
\ET

\BR\label{rm-Primal}
Alternatively, we can estimate the growth of function $\xi_{\hat z_k}(\cdot)$ using the second representation in (\ref{eq-OGrow1}).
For that, we need to prove that $\| \Delta x_k \|_{\hat x_k} \leq O(\sqrt{\nu_g})$ and repeat the reasoning before Theorem \ref{th-DBS1}. We move these justifications in Appendix in Section \ref{sc-App} since they are more technical due to a non-trivial way of generating the point $\hat x_k$. This {\em primal} justification gives better results if $\nu_g \leq \nu_{\zeta}$. However, in both cases the method (\ref{met-DBS}) is not changing. We just use different formulas for computing values of the same primal-dual proximity measure.
Therefore, in the sequel we claim that the complexity of method (\ref{met-DBS}) for symmetric cones depends on $\min\{ \nu_{\zeta}, \nu_g \}$. \qed
\ER

Thus, method (\ref{met-DBS}), being applied to a conic problem (\ref{prob-PD}) with symmetric cone $K$, solves it in $O\left(\min\{\nu_{\zeta}, \nu_g \}^{1/2}\ln {1 \over \epsilon}\right)$ iterations of the dual Newton method. In Section \ref{sc-Easy}, we will discuss some classes of problems with $\nu_{\zeta} << \nu$. 

To conclude this section, let us present a variant of method (\ref{met-DBS}), where we eliminate all unnecessary appearances of the primal variables.

\beq\label{met-DualBS}
\ba{|l|}
\hline \\
\hspace{5ex} \mbox{\sc Dual Predictor-Corrector Method for Symmetric Cones}\\
\\
\hline \\
\mbox{{\bf Initialization.} Let $\beta$ and $\delta$ satisfy (\ref{eq-condDBS1}). Define $A = \delta + \omega_*\left({2 \beta \over (1-\beta)^2}\right)$. }\\
\\
\mbox{Compute $y_0 \in {\cal F}_y$ with $\| \nabla \zeta(y_0) \|_{y_0} \leq \half \beta$ and set $t_0 = {\beta - \| \nabla \zeta(y_0) \|_{y_0} \over  \| b \|_{y_0}}$.}\\
\\
\mbox{\bf $k$-th iteration ($k \geq 0$)}\\ \\
\ba{rl}
\mbox{a)} & \mbox{For $y_k \in \inter {\cal F}_y$, compute $g_k = \nabla \psi_{t_k}(y_k)$ and $B_k=[ \nabla^2 \zeta(y_k)]^{-1}$.}\\
\\
& \mbox{Define  $d_k = B_k g_k$, and $\lambda_k = \la g_k, d_k \ra^{1/2}$.}\\
\\
\mbox{b)} & \mbox{{\bf If} $\lambda_k > \beta$, {\bf then} $y_{k+1} = y_k - {d_k \over 1 + \lambda_k}$, $t_{k+1} = t_k$ ({\em Damped Newton Step})}\\
\\
\mbox{c)} & \mbox{{\bf else} ({\em Predictor Step})}\\ \\
& \mbox{1. Set $s_k = c - A^* y_k$, $\hat y_k = y_k + d_k$, $\hat s_k = c - A^* \hat y_k$, and $\Delta y_k = t_k B_k b$.}\\
\\
& \mbox{2. Define function \fbox{$\xi_k(\alpha) = \zeta(\hat y_k + \alpha \Delta y_k) + \zeta\left(\hat y_k - {\alpha \over 1 - \alpha} \Delta y_k \right) - 2 \zeta(\hat y_k)$}}\\
\\
& \mbox{3. Compute the step $\alpha_k \in (0,1)$ satisfying the equation $\xi_k(\alpha_k) = A$.}\\
\\
& \mbox{4. Define $y_{k+1} = \hat y_k + \alpha_k \Delta y_k$ and $t_{k+1} =  {\nu \, \cdot \, t_k \over (1-\alpha_k)  \| \hat s_k \|^2_{s_k}}$.}
\ea\\
\\
\hline
\ea
\eeq

In this scheme, as compared with method (\ref{met-DBS}), the only comments are needed for the Step~c$_4$). Indeed,
$$
\ba{rcl}
\la \hat s_k, \hat x_k \ra & \refEQ{eq-Out} & \la \hat s_k, \nabla^2 F_*(w^*_k) \hat s_k \ra \; \refEQ{met-DBS} \; {1 \over t_k} \la \hat s_k, \nabla^2 F_*(s_k) \hat s_k \ra \; = \; {1 \over t_k} \| \hat s_k \|^2_{s_k}.
\ea
$$ 
The last equality in this chain can be used at Step c4) for computing $t_{k+1}$.
For termination criterion, we can verify the condition
\beq\label{eq-Term}
\ba{rcl}
{\nu \over t_{k+1}} & \leq & \epsilon
\ea
\eeq
in the end of the predictor step.
If it is fulfilled, we generate the output $\hat z_k + \alpha_k \Delta z_k \in {\cal F}$  by the rules of algorithm (\ref{met-DBS}).

\section{Easy Linear Matrix Inequalities}\label{sc-Easy}
\SetEQ

Let us demonstrate the main advantages of method (\ref{met-DualBS}) as applied to problems of Semidefinite Optimization. We start from the following simple result.
\BL\label{lm-NuSDP}
Let matrix $C \in \R^{n \times n}$ be positive definite and matrix $A \in \R^{m \times n}$, $m < n$, has full row rank. 
Then the self-concordant barrier 
\beq\label{def-ZY}
\ba{rcl}
\zeta(y) & \Def & - \ln \det (C - A^T D(y) A), \quad D(y)  = \Diag(y),
\ea
\eeq
has the following short representation
\beq\label{eq-RepZ}
\ba{rcl}
\zeta(y) & = & - \ln \det(G^{-1} - D(y)) - \ln (\sigma \gamma),
\ea
\eeq 
where $G = A C^{-1} A^T$,
$\sigma = \det C$, and $\gamma = \det G$. For its domain
$$
\ba{rcl}
\dom \zeta(y) & = & \{ y \in \R^m: \; D(y) \prec G^{-1} \},
\ea
$$
the parameter
$\nu_{\zeta}$ has the best possible value:
\beq\label{eq-NuBound}
\ba{rcl}
\nu_{\zeta} & = & m.
\ea
\eeq
\EL
\proof
Indeed, for $U \Def G^{-1/2} A C^{-1/2}$, we have $U U^T = I_m$. Therefore
$$
\ba{rcl}
\det(C - A^T D(y) A) & = & \sigma \det(I_n - C^{-1/2} A^T D(y) A C^{-1/2}) \\
\\
& = &
\sigma \det (I_n - U^T G^{1/2} D(y) G^{1/2} U ) \; = \;
\sigma  \det (I_m - G^{1/2} D(y) G^{1/2}) \\
\\
& = & \sigma \gamma \det[ G^{-1} - D(y)]. 
\ea
$$

Since the barrier parameter of function $- \ln \det Y$, $Y = Y^T \in \R^{m \times m}$, is equal to $m$, we get $\nu_{\zeta} \leq m$. Let us prove now the inverse inequality.

Denote ${\cal F}_y = \{ y \in \R^m: \; G^{-1} - D(y) \succeq 0\}$.
For $\tau > 0$, let us choose $\bar y = - \tau e$. Note that $p_i \Def -e_i$, $i=1,\dots, m$, are recession directions of the feasible set ${\cal F}$. For each $i$, define the value $\beta_i$ from the condition
$$
\ba{rcl}
\bar y - \beta_i p_i \; = \; - \tau e + \beta_i e_i & \not\in & \inter {\cal F}.
\ea
$$
That is $G^{-1} + \tau I_m \not\succ \beta_i e_i e_i^T$, which gives us $\beta^{-1}_i = \la (G^{-1}+\tau I_m)^{-1} e_i, e_i \ra$.
On the other hand, for $\alpha_i =  \alpha \Def\lambda_{\min}(G^{-1}+\tau I_m)$, we have
$$
\ba{rcl}
G^{-1} + \tau  I_m & \succeq & \alpha I_m \; = \; \sum\limits_{i=1}^m \alpha_i e_i e_i^T.
\ea
$$
This means that $\bar y - \sum\limits_{i=1}^m \alpha_i p_i = - \tau  e + \alpha e \in {\cal F}_y$. Therefore, in view of Theorem \ref{th-LowB}, we have
$$
\ba{rcl}
\nu_{\zeta} & \geq & \sum\limits_{i=1}^m {\alpha_i \over \beta_i} \; = \; \lambda_{\min}(G^{-1} + \tau I_m) \la (G^{-1} + \tau I_m)^{-1}, I_m \ra \\
\\
& = & \lambda_{\min}(\tau^{-1}G^{-1} + I_m) \la (\tau^{-1}G^{-1} + I_m)^{-1}, I_m \ra.
\ea
$$ 
Taking the limit in the right-hand side of this inequality as $\tau \to + \infty$, we get $\nu_{\zeta} \geq m$.
\qed

Note that the barrier parameter of function $\zeta(\cdot)$ does not depend on its representation. Thus, any of representations (\ref{def-ZY}) or (\ref{eq-RepZ}) can be used for numerical computations. The short form (\ref{eq-RepZ}) is computationally more preferable, but it needs some pre-processing.

The result of Lemma \ref{lm-NuSDP} can be easily extended in the following way.
\BC\label{cor-Par}
Let matrix $C \in \R^{n \times n}$ be positive definite and the matrices $A_i \in \R^{n \times n}$, $i = 1, \dots, m$, be symmetric. Then the barrier parameter of self-concordant function
$$
\ba{rcl}
\zeta(y) & = & - \ln \det \Big( C - \sum\limits_{i=1}^m y^{(i)} A_i \Big)
\ea
$$
satisfies the following bound:
\beq\label{eq-BGen}
\ba{rcl}
\nu_{\zeta} & \leq & \min \Big\{ n, \sum\limits_{i=1}^m \mbox{\rm rank}(A_i) \Big\}.
\ea
\eeq
\EC

Clearly, the right-hand side of inequality (\ref{eq-BGen}) is minimal when the rank of all matrices $A_i$ is equal to one. In any case, this
result is useful for treating problems of Semidefinite Optimization (\ref{prob-SDOP}), (\ref{prob-SDOD})
with {\em low-rank} equality constraints, especially when $n \to \infty$. It opens an interesting possibility for constructing polynomial-time methods for some problems of {\em Semi-Infinite Optimization}.

It appears that for Semidefinite Optimization, the dual representation of the growth of functional proximity measure (\ref{eq-OGrow}) admits further simplification. Indeed, in the case
$$
\ba{rcl}
F_*(S) & = & - \ln \det S - \nu,
\ea
$$
we can represent $S = L_n L_n^T \in \R^{n \times n}$, where $L$ is a lower-triangular matrix.  Denote 
$$
\ba{rcl}
B & = & L_n^{-1} \Delta S L_n^{-T}.
\ea
$$
Then
\beq\label{eq-RepB}
\ba{rcl}
\xi_z(\alpha) & = & - \ln \det(S + \alpha \Delta S) - \ln \det \left(S - {\alpha \over 1 - \alpha} \Delta S\right) + 2 \ln \det S\\
\\
& = & - \ln \det(I_n + \alpha B) - \ln \det \left(I_n - {\alpha \over 1 - \alpha} B\right) \\
\\
& = & - \ln \det \left(I_n -  {\alpha^2 \over 1 - \alpha}(B + B^2) \right).
\ea
\eeq
Further, in $O(n^3)$ we can compute representation $B=U_n T_n U_n^T$, where $U_n$ is a unitary matrix and matrix $T_n \in \R^{n \times n}$ is three-diagonal. Then
\beq\label{eq-RepT}
\ba{rcl}
\xi_z(\alpha) & = & - \ln \det(I_n + \alpha T_n) - \ln \det \left(I_n - {\alpha \over 1 - \alpha} T_n\right) \\
\\
& = & - \ln \det \left(I_n -  {\alpha^2 \over 1 - \alpha}(T_n + T_n^2) \right).
\ea
\eeq
Note that matrix $T_n + T_n^2$ has five nonzero diagonals. Therefore, one computation of the value of function $\xi_z(\cdot)$ by any representation in (\ref{eq-RepT}) needs $O(n)$ arithmetic operations. Hence, the complexity of solving equation at Step c3) of method (\ref{met-DualBS}) by a bisection scheme is $O(n\ln\epsilon_a^{-1})$ arithmetic operations, where $\epsilon_a$ is the requried accuracy for satisfying the target equation.

A similar trick works for the low-rank constraints. Assume for simplicity that all equality constraints in the problem (\ref{prob-SDOP}$_1$) have rank one:
\beq\label{eq-ROne}
\ba{rcl}
A_i & = & a_i a_i^T \in \R^{n \times n}, \quad i = 1, \dots, m.
\ea
\eeq
Then, in accordance to Lemma \ref{lm-NuSDP}, we have the following representation of function $\xi_z(\cdot)$:
$$
\ba{rcl}
\xi_z(\alpha) & = & - \ln \det (G^{-1} - \alpha \Delta Y) - \ln \det \left(G^{-1} + { \alpha \over 1 - \alpha} \Delta Y \right) + 2 \ln \det G^{-1},
\ea
$$
where  $\Delta Y = \Diag (\Delta y) \in \R^{m \times m}$. Denoting $B_m = L^T_m \Delta Y L_m$, where $ L_m L_m^T=G$, we can rewrite it as follows:
$$
\ba{rcl}
\xi_z(\alpha) & = & - \ln \det(I_m + \alpha B_m) - \ln \det \left(I_m - {\alpha \over 1 - \alpha} B_m\right).
\ea
$$ 
It remains to compute factorization $B_m = U_m T_m U_m^T$ with a unitary matrix $U_m$ and three-diagonal matrix $T_m$. Then, we get the following representations:
\beq\label{eq-RepTM} 
\ba{rcl}
\xi_z(\alpha) &  = & - \ln \det(I_m + \alpha T_m) - \ln \det \left(I_m - {\alpha \over 1 - \alpha} T_m\right)\\
\\
& = & - \ln \det \left(I_m + {\alpha^2 \over 1 - \alpha} (T_m - T_m^2) \right).
\ea
\eeq
Computation of the value of the right-hand side in the representations (\ref{eq-RepTM}) needs $O(m)$ operations.

Let us show that for SDO-problem (\ref{prob-SDOP}) with low-rank constraints, the dual problem is indeed much simpler.
\BE\label{ex-LowRank}
Consider the following primal-dual SDO-problem:
\beq\label{prob-ExLR}
\ba{rcl}
\min\limits_{X \succeq 0} \{ \la I_n, X \ra: \; \la X a, a \ra = 1 \} & = & \max\limits_{y \in \R} \{ y: \; I_n \succeq y \cdot aa^T \},
\ea
\eeq
where $X = X^T \in \R^{n \times n}$ and $\| a \| = 1$. Note that the feasible set of the dual problem is unbounded. Hence, we cannot rely on Lemma \ref{lm-FY}. However, in accordance to Lemma \ref{lm-NuSDP}, for solving the dual problem we can use the following barrier function:
$$
\ba{rcl}
\zeta(y) & = & - \ln \det (I_n - y a a^T) \; = \; - \ln(1-y)
\ea
$$
with $\nu_{\zeta} = 1$. At the same time, we can prove that the value of parameter $\nu_g$ of any self-concordant barrier $g(\cdot)$ for the primal feasible set
$$
\ba{rcl}
{\cal F}_p & = & \{ X \succeq 0: \; \la X a, a \ra = 1 \}
\ea
$$
cannot be less than $n-1$. Indeed, without loss of generality, we can assume that $a = e_1$. Let us choose $\bar X \Def I_n \in {\cal F}_p$.
It is clear that $p_i = e_i e_i^T$, $i = 2, \dots , n$, are recession directions of the set ${\cal F}_p$. At the same time, for $\beta_i = 1$, $i = 2, \dots, n$, we have
$$
\ba{rcl}
\bar X - \beta_i p_i & = & I_n - e_i e_i^T \; \not\in \; \inter {\cal F}_p.
\ea
$$
On the other hand, for $a_i = 1$, $i=2, \dots ,n$, we have
$$
\ba{rcl}
\bar X - \sum\limits_{i=2}^n a_i p_i & = & e_1 e_1^T \; \in {\cal F}_p.
\ea
$$
Therefore, in view of Theorem \ref{th-LowB}, we get
$\nu_g \geq \sum\limits_{i=2}^n {a_i \over b_i} \; = \; n-1$. \qed
\EE

A similar effect can also occur in some problems of Linear Optimization.
\BE\label{ex-Lin}
Let us choose dimension $n \geq 1$ and vector $b \in \R^n$ with coefficients
$$
\ba{rcl}
b^{(1)} & = & 1, \quad b^{(i)} \; \in \; (0,1), \quad i = 2, \dots, n.
\ea
$$
Denote 
$a \Def (b^T,1,-b^T)^T \in \R^{2n+1}$. Consider the following pair of primal-dual problems:
\beq\label{prom-ELin}
\ba{rcl}
\min\limits_{x \in \R^{2n+1}_+} \{ \la e, x \ra: \; \la a, x \ra = 1 \} & = & \max\limits_{y \in \R} \{ y: \; e \geq y \cdot a \}.
\ea
\eeq
For the primal feasible set ${\cal F}_p = \{ x \geq 0: \la a, x \ra = 1\}$, we have the following recession directions: 
$$
\ba{rcl}
p_i & = & (e_i^T, 0, e_i^T)^T, \quad i = 1, \dots , n.
\ea
$$
At the point $\bar x \Def e \in \R^{2n+1}$, we can see that
$$
\ba{rcl}
\bar x - \sum\limits_{i=1}^n p_i & = & e_{n+1} \; \in \; {\cal F}_p, \quad \bar x - p_i \; \not\in \; \inter {\cal F}_p, \quad i = 1, \dots, n.
\ea
$$
Hence, by Theorem \ref{th-LowB}, any self-concordant barrier $g(\cdot)$ for set ${\cal F}_p$ has the barrier parameter $\nu_g \geq n$. This bound does not depend on the coefficients of vector $b$.

At the same time, the standard logarithmic barrier of the primal cone $K = \R^{2n+1}_+$ is
$$
\ba{rcl}
F(x) & = & - \sum\limits_{i=1}^{2n+1} \ln x^{(i)}, \quad \nu_F \; = \; 2 n+1,
\ea
$$
with the dual barrier $F_*(s) = - \sum\limits_{i=1}^{2n+1} \ln s^{(i)} - \nu_F$. Hence, the barrier function for the dual feasible set ${\cal F}_d \Def (-1,1)$ is
$$
\ba{c}
\zeta(y) \; = \; F_*(e - y \cdot a ) \; = \; \zeta_1(y) + \zeta_2(y),\\
\\
\zeta_1(y) \; \Def \;
- \ln(1 - y) - \ln(1-y^2),  \quad
\zeta_2(y) \; \Def \; - \sum\limits_{i=2}^n \ln \left(1 - (b^{(i)} y)^2 \right).
\ea
$$
Note that $\zeta_1(\cdot)$ is a self-concordant barrier for ${\cal F}_d$ with parameter $\nu_1 = 2$ and 
$$
\ba{rcl}
\zeta_1''(y) & = & {2 \over (1-y)^2} + {1 \over (1+y)^2} \; \geq \; 2, \quad y \in (-1,1).
\ea
$$
For the second barrier function, we have
$$
\ba{rcl}
|\zeta_2'(y)| & = & 2 \sum\limits_{i=2}^n {(b^{(i)})^2 |y| \over 1 - (b^{(i)}y)^2}  \; \leq \; B \; \Def \; 2 \sum\limits_{i=2}^n {(b^{(i)})^2 \over 1 - (b^{(i)})^2}, \quad y \in {\cal F}_d.
\ea
$$
Hence,
$$
\ba{rcl}
\nu_{\zeta} & = & \sup\limits_{|y|<1} {(\zeta'(y))^2 \over \zeta''(y)} \; \leq \; 2 \sup\limits_{|y|<1} {(\zeta_1'(y))^2 + (\zeta_2'(y))^2 \over \zeta_1''(y) + \zeta''_2(y)} \; \leq \; 2 \left( 2 + \sup\limits_{|y|<1} {(\zeta_2'(y))^2 \over \zeta_1''(y) } \right)\; \leq \; 2 \left(2+ \half B^2 \right).
\ea
$$
Thus, choosing the coefficients of vector $b$ small enough, we can have $B \to 0$. Consequently, the barrier parameter $\nu_{\zeta}$ will be as close to four as we wish. Again, this is a situation when the parameters of self-concordant barriers for the primal and the dual feasible sets are fundamentally different.
\qed
\EE

To conclude this section, let us show that for functional proximity measures of other symmetric cones there also exist short representations similar to (\ref{eq-RepB}).

{\bf 1.} For $K = K^* = \R^n_+$, we choose $F_*(s) = -\sum\limits_{i=1}^n \ln s^{(i)}+n$. Then the definition (\ref{def-TanD}) is as follows:
\beq\label{def-ALP}
\ba{rcl}
{\Delta s \over s} + {\Delta x \over x} & = & - e, \quad A \Delta x = 0, \quad \Delta s + A^T \Delta y = 0.
\ea
\eeq
Therefore,
\beq\label{eq-ShortLP}
\ba{rcl}
\xi_z(\alpha) & \refEQI{eq-OGrow}{2} & - \sum\limits_{i=1}^m \Big[\ln (s^{(i)} + \alpha \Delta s^{(i)}) + \ln \left(s^{(i)} - {\alpha \over 1 - \alpha} \Delta s^{(i)} \right) - 2 \ln s^{(i)} \Big]\\
\\
& = & - \sum\limits_{i=1}^m \ln \Big[ \left(1 + \alpha {\Delta s^{(i)} \over s^{(i)}}\right) \cdot \left( 1- {\alpha \over 1 - \alpha} \cdot {\Delta s^{(i)} \over s^{(i)}}\right) \Big]\\
\\
& = & - \sum\limits_{i=1}^m \ln \left( 1- {\alpha^2 \over 1 - \alpha}{\Delta s^{(i)} \over s^{(i)}} - {\alpha^2 \over 1 - \alpha}
\left({\Delta s^{(i)} \over s^{(i)}}\right)^2 \right) \\
\\
& \refEQ{def-ALP} & - \sum\limits_{i=1}^m \ln \left( 1 + {\alpha^2 \over 1 - \alpha}{\Delta s^{(i)} \over s^{(i)}}{\Delta x^{(i)} \over x^{(i)}}\right).
\ea
\eeq

{\bf 2.} For the Lorentz cone $K^* = \{ s = (\tau,u) \in \R_+ \times \R^{n}: \; \tau^2 \geq \| u \|^2 \}$, we have 
$$
\ba{rcl}
F_*(s) = - \ln (\tau^2 - \| u \|^2) - 2 + 2 \ln 2.
\ea
$$
Therefore, denoting $\omega = \tau^2 - \| x \|^2$, we get
$$
\ba{rcl}
\xi_z(\alpha) & \refEQI{eq-OGrow}{2} & - \ln (\omega + 2\alpha \Delta_1 + \alpha^2 \Delta_2) - \ln \left(\omega - {2\alpha \over 1 - \alpha} \Delta_1 + {\alpha^2 \over (1- \alpha)^2} \Delta_2\right) + 2 \ln \omega\\
\\
& = & - \ln (1 + 2\alpha \delta_1 + \alpha^2 \delta_2) - \ln \left(1 - {2\alpha \over 1 - \alpha} \delta_1 + {\alpha^2 \over (1- \alpha)^2} \delta_2\right),
\ea
$$
where $\Delta_1 = \tau \Delta \tau - \la u, \Delta u \ra$, $\delta_2 = \Delta \tau^2 - \| \Delta u \|^2$, and $\delta_i = {1 \over \omega} \Delta_i$, $i=1,2$. Note that
$$
\ba{rl}
& (1 + 2\alpha \delta_1 + \alpha^2 \delta_2)\cdot \left(1 - {2\alpha \over 1 - \alpha} \delta_1 + {\alpha^2 \over (1- \alpha)^2} \delta_2\right)\\
\\
= & 1 + 2\alpha \delta_1 + \alpha^2 \delta_2 - {2\alpha \over 1 - \alpha} \delta_1 -  {4\alpha^2 \over 1 - \alpha} \delta^2_1 - {2\alpha^3 \over 1 - \alpha} \delta_1 \delta_2 + {\alpha^2 \over (1- \alpha)^2} \delta_2 + {2\alpha^3 \over (1- \alpha)^2} \delta_1 \delta_2 + {\alpha^4 \over (1- \alpha)^2} \delta^2_2\\
\\
= & 1  - {2\alpha^2 \over 1 - \alpha} \delta_1 + \alpha^2 \left(1 + {1 \over (1- \alpha)^2} \right)\delta_2  + {2\alpha^3 \over 1- \alpha}\left( {1 \over 1 - \alpha} - 1 \right) \delta_1 \delta_2 -  {4\alpha^2 \over 1 - \alpha} \delta^2_1  + {\alpha^4 \over (1- \alpha)^2} \delta^2_2\\
\\
= & 1  - {2\alpha^2 \over 1 - \alpha} [\delta_1 + 2 \delta_1^2 - \delta_2]+  {\alpha^4 \over (1- \alpha)^2} [\delta^2_2 + 2 \delta_1 \delta_2 + \delta_2] .
\ea
$$
Thus, for $a_1 \Def \delta_1 + 2 \delta_1^2 - \delta_2$ and $a_2 \Def \delta^2_2 + 2 \delta_1 \delta_2 + \delta_2$, we get the following representation:
\beq\label{eq-ShortLR}
\ba{rcl}
\xi_z(\alpha) & = & - \ln\left(1 - {2\alpha^2 \over 1 - \alpha} a_1 + \left( {\alpha^2 \over 1 - \alpha} \right)^2 a_2 \right). 
\ea
\eeq

As we have seen in our examples (\ref{eq-RepB}), (\ref{eq-ShortLP}), and (\ref{eq-ShortLR}), the most convenient search variable for ensuring condition $\xi_k(\alpha_k) = A$ at Step c$_3$) of the method (\ref{met-DualBS}) is $\tau = {\alpha^2 \over 1 - \alpha}$. We can expect that in the neighborhood of optimal solution, the coefficients for variable $\tau$ in the above representations are vanishing. This opens a possibility for $\alpha_k \to 1$, which results in the local superlinear convergence of the method.

\section{Numerical experiments}\label{sc-Num}
\SetEQ

For numerical experiments, we select a nonsmooth optimization problem in the space of symmetric matrices, which we call the {\em Low-Rank Quadratic Interpolation}:
\beq\label{prob-LRQI}
\ba{c}
\min\limits_{X = X^T \in \R^{n \times n}} \Big\{ \sum\limits_{i=1}^n |\lambda_i(X)|: \; \la X a_i, a_i \ra = b_i, \quad i = 1, \dots, m \Big\}.
\ea
\eeq
In this formulation, the right-hand sides of the equality constraints can be of any sign. We assume that the matrix $A^T \Def (a_1, \dots, a_m)$ has full column rank. Thus, $G \Def AA^T \succ 0$.

This problem admits the following SDO-formulation:
\beq\label{prob-IntSDP}
\ba{c}
\min\limits_{X_1, X_2 \succeq 0} \Big\{ \la I_n, X_1 + X_2 \ra :\; 
\la (X_1 - X_2) a_i, a_i \ra = b_i, \quad i = 1, \dots, m \Big\}.
\ea
\eeq
The dual problem to (\ref{prob-IntSDP}) is as follows:
\beq\label{prob-IntDual}
\ba{rcl}
f_* & = & \max\limits_{y \in \R^m} \Big\{ \la b, y \ra: \; - I_n \preceq A^T D(y) A \preceq I_n \Big\},
\ea
\eeq
where  $D(y) = \Diag(y) \in \R^{m \times m}$. Thus, we can choose the following dual barrier function
\beq\label{eq-DualF}
\ba{rcl}
\zeta(y) & = & - \ln \det (I_n - A^T D(y) A) - \ln \det (I_n + A^T D(y) A)\\
\\
& \refEQ{eq-RepZ} & - \ln \det (G^{-1} - D(y))  - \ln \det (G^{-1} + D(y)) - 2 \ln \sigma,
\ea
\eeq
where $\sigma = \det G$. Since the feasible set of problem (\ref{prob-IntDual}) is bounded, the Hessian $\nabla^2 \zeta(\cdot)$ is positive definite at any feasible point.
Note that 
$$
\ba{rcl}
\zeta(0) & = & 0, \quad \nabla \zeta(0) \; = \; 0.
\ea
$$
Therefore, for method (\ref{met-DualBS}) we choose $y_0 = 0$.

In accordance to Lemma \ref{lm-FY}, 
the parameter of barrier function for the primal feasible set is big: $\nu_g = 2n$. 
Let us show that $\nu_{\zeta} = m$. This is a consequence of the following result.
\BL\label{lm-MBox}
The barrier $F(Y) = - \ln \det (I_m -Y) - \ln \det(I_m +Y)$ for the matrix box 
\beq\label{eq-MBox}
\ba{rcl}
{\cal B}_m & = & \{ Y = Y^T \in \R^{m \times m}: \; - I_m \preceq Y \preceq I_m \}
\ea
\eeq
has parameter $\nu_{F} = m$.
\EL
\proof
Let us compute the first derivative of $F(\cdot)$ at some feasible $Y$ along direction $H$, which is a symmetric $m \times m$-matrix:
$$
\ba{rcl}
\la \nabla F(Y), H \ra & = & \la (I_m-Y)^{-1} - (I_m+Y)^{-1}, H \ra\\
\\
& = & 2 \la (I_m-Y)^{-1} Y (I_m+Y)^{-1}, H \ra \; = \;
\la D\left({2 \lambda \over 1 - \lambda^2} \right), \hat H \ra,
\ea
$$
where we represent $Y = U D(\lambda) U^T$ by eigenvalue decomposition and denote $\hat H = U^T H U$. Similarly, for the second derivative we have
$$
\ba{rcl}
\la \nabla^2 F(Y) H, H \ra & = & \la (I_m-Y)^{-1}H(I_m-Y)^{-1} + (I_m+Y)^{-1}H (I_m+Y)^{-1}, H \ra\\
\\
& = & \la D \left({1 \over 1- \lambda}\right) \hat H D \left({1 \over 1- \lambda}\right), \hat H \ra + \la D \left({1 \over 1+ \lambda}\right) \hat H D \left({1 \over 1+ \lambda}\right), \hat H \ra\\
\\
& = & \la \tilde H {1 \over 1-\lambda}, {1 \over 1-\lambda} \ra + \la \tilde H {1 \over 1+\lambda}, {1 \over 1+\lambda} \ra\; \geq \; \la d^2,  {2(1 + \lambda^2) \over (1-\lambda^2)^2}  \ra \; \geq \; \la d^2, \left({2 \lambda \over 1 - \lambda^2} \right)^2 \ra,
\ea
$$
where the matrix $\tilde H$ has elements $\tilde H_{i,j} = \hat H_{i,j}^2 \geq 0$ and $d = \diag(\hat H)$. Hence, by Cauchy-Schwartz inequality, we have
$$
\ba{rcl}
\la \nabla^2 F(Y) H, H \ra & \geq & \la d^2, \left({2 \lambda \over 1 - \lambda^2} \right)^2 \ra \; \geq {1 \over m} \la |d|, {2 \lambda \over 1 - \lambda^2} \ra^2 \; \geq \; {1 \over m} \la \nabla F(Y), H \ra^2.
\ea
$$
This means that $\nu_F \refEQ{eq-SCB} m$.
\qed

We tested method (\ref{met-DualBS}) on several series of one hundred random problems (\ref{prob-LRQI}). In our random generator, we choose $m$ and $n\geq 2m$, and generate random vectors $a_i \in \R^n$, $i=1 \dots , m$, and vector $b \in \R^m$ with entries uniformly distributed in the interval $[-1,2]$. For parameters of the method, we choose
$$
\ba{rcl}
\beta & = & 0.2, \quad A \; = \; 2, \quad \epsilon \; = \; 10^{-8}, \quad y_0 \; = \; 0 \; \in \; \R^m.
\ea
$$
The parameter of self-scaled barrier for the primal cone in problem (\ref{prob-IntSDP}) is $\nu = 2n$. We use it at Step c$_4$) of the method (\ref{met-DualBS}).

Below, in tables (\ref{tab-ResP}), (\ref{tab-ResP}), and (\ref{tab-Super}), we present our computational results. In the first two tables, each line corresponds to different value of $m$, and each column corresponds to different $n$. In each cell, we show an average number of iterations with relative standard deviation.  

At this moment, the most interesting results are presented in table (\ref{tab-ResP}), where we show the number of predictor steps. These numbers, despite to the high accuracy $\epsilon = 10^{-8}$, are very small (from six to ten). This means that the functional proximity measure, employed at the Step c$_2$) of method (\ref{met-DualBS}), indeed admits very large steps, adjusting to the actual curvature of the central path.
\beq\label{tab-ResP}
\ba{c}
\mbox{Average number of Predictor Steps (100 problems)}\\
\ba{|c|c|c|c|c|c|}
\hline
m \; \backslash \; n & 64 & 128 & 256 & 512 & 1024 \\
\hline
32 & 9.0 \pm 9.6 \% & 8.2 \pm 9.3 \% & 7.1 \pm 4.2 & 6.9 \pm 3.5\%  & 6.6 \pm 7.4 \% \\ \hline 
64 & & 9.9 \pm 7.8\% & 7.8 \pm 5.5\% &  7.1 \pm 3.4\% & 6.9 \pm 3.9\% \\ \hline
128 &  & & 9.9 \pm  6.2\% & 7.9 \pm 4.3\% & 7.0 \pm 2.0\% \\ \hline
256 & & & & 9.5 \pm 5.7\% & 7.9 \pm 3.7\% \\ \hline
512 & & & & & 9.5 \pm 5.3 \% \\
\hline
\ea
\ea
\eeq

In the second table (\ref{tab-ResK}), we show the average values for the total number of steps of the method. As we can see, each predictor step needs approximately four corrector steps in order to come back to the neighborhood of the central path. This number is not too big and it is clear that it can be reduced by more sophisticated centering strategies. In our experiments, we just used the rudimentary Damped Newton Method. Its improvement is an interesting topic for further research. 
\beq\label{tab-ResK}
\ba{c}
\mbox{Average number of iterations (100 problems)}\\
\ba{|c|c|c|c|c|c|}
\hline
m \; \backslash \; n & 64 & 128 & 256 & 512 & 1024 \\
\hline
32 & 40.9 \pm 13.6 \% & 37.2 \pm 13.9 \% & 31.8 \pm 7.3 & 32.0 \pm 6.0\%  & 30.5 \pm 7.4 \% \\ \hline 
64 & & 47.2 \pm 10.3\% & 35.0 \pm 5.9\% &  32.3 \pm 4.8\% & 32.9 \pm 4.7\% \\ \hline
128 &  & & 48.3 \pm  8.8\% & 35.5 \pm 5.4\% & 33.0 \pm 3.5\% \\ \hline
256 & & & & 46.3 \pm 7.1\% & 36.2 \pm 5.0\% \\ \hline
512 & & & & & 46.0 \pm 6.2 \% \\
\hline
\ea
\ea
\eeq

In the last table (\ref{tab-Super}), we show the typical dynamics for convergence of method~(\ref{met-DualBS}). Usually, in the end of the process, the rate of convergence becomes superlinear. At the same time, it needs more and more cheap bisection steps in the line search procedure of Step c$_3$). Note that, after an appropriate tri-duaginalization of $m \times m$-matices, the cost of such a step is only $O(m)$.
\beq\label{tab-Super}
\ba{c}
\mbox{Progress by iterations ($m=256$, $n=1024$)}\\
\ba{|c|c|c|c|}
\hline
N_{\rm pred}/N_{\rm tot} & \la \hat s_k, \hat x_k \ra & t_k & \mbox{Bisections} \\
\hline 
0/0 &  & 31.25 &  \\
\hline
1/2 & 9.7  \cdot 10^1 & 2.4 \cdot 10^{2} & 9 \\
\hline
2/6 & 8.6 \cdot 10^0 & 7.1 \cdot 10^{2} & 7 \\
\hline
3/9 & 2.9 \cdot 10^0 & 2.5 \cdot 10^{3} & 8\\
\hline
4/14 & 8.4 \cdot 10^{-1} & 1.3 \cdot 10^{4} &  8 \\
\hline
5/21 & 1.5 \cdot 10^{-1} & 2.6 \cdot 10^{5} & 10  \\
\hline
6/27 & 7.9 \cdot 10^{-3} & 6.8 \cdot 10^{7} & 13 \\
\hline
7/33 & 3.0 \cdot 10^{-5} & 8.9 \cdot 10^{10} & 15\\
\hline
8/38 & 6.8 \cdot 10^{-9} & 3.2 \cdot 10^{11} &  \\
\hline
\ea
\ea
\eeq

In our opinion, the above computational results are very promising. They show that some SDO-problems can be solved in the time comparable with the time required by Linear Optimization. We performed our experiments at the usual notebook. Nevertheless, the largest problem in our test set was solved in less than fifty seconds. For implementing all expensive operations of Linear Algebra (inverting matrices, computing determinant, etc.) we used the standard Cholesky factorization.

It will be interesting to compare our computational results with  performance of other efficient optimization schemes (e.g. \cite{LST,ZST}). For that, we present a reformulation of problem~(\ref{prob-IntDual}), suitable for solving by the other first- and second-order methods:
\beq\label{prob-IntRef}
\ba{rcl}
f_*^{-1} & = & \min\limits_{y \in \R^m} \Big\{ \; \sigma_{\max} \left( \sum\limits_{i=1}^m y^{(i)} a_i a_i^T \right): \; \la b, y \ra = 1 \; \Big\}.
\ea
\eeq

\section{Conclusion}\label{sc-Conc}
\SetEQ

Let us summarize the main results of the paper.

{\bf 1.} We proposed a new methodology for solving primal-dual problems of Conic Optimization by predictor-corrector schemes. Instead of finding at the Corrector Stage a primal-dual pair in a close neighborhood of the primal-dual central path, we run the corrector process only in the dual space. We use the obtained approximation of the dual path as a scaling point for a new feasible primal-dual pair, which can be easily generated by the Dual Gambit Rule (\ref{eq-Out}).
After that, we perform a predictor step using the standard affine-scaling direction with the step size defined by {\em Functional Proximity Measure}. The computational complexity of this technique is very attractive even for SDO-problems, where we need to apply only the standard Cholesky factorizations.

{\bf 2.} One of the main motivations for our approach is the observation that in many situations the complexity of primal and dual problems are different. For symmetric cones, our method is the first one with the complexity bounds dependent on the minimal value between the barrier parameters for the primal and the dual feasible sets. We prove that the standard assumption on boundedness of the dual feasible set necessarily results in the maximal value of barrier parameter for the primal feasible set. At the same time, as we show by examples, the barrier parameter for the dual feasible set can be much smaller. The automatic switching to the best value of the barrier parameter is a very desirable property. In our case, it is ensured by the use of functional proximity measure.

{\bf 3.} We discuss applications of this technique to several problem of Semidefinite Optimization. It appears that computational complexity of some of them is similar to that of the problems of Linear Optimization with the comparable dimension. 

{\bf 4.} For our approach, we present the results of preliminary testing on randomly generated problems of Low-Rank Quadratic Interpolation. They are quite promising since the number of predictor steps is small and we observe in the end a local quadratic convergence.

It will be very interesting to check the performance of the proposed approach on other structural problems of Conic Optimization. One of the most important open questions is the possibility of having asymmetric complexity bounds for other cones and other primal-dual proximity measures. The second part of this question is open even for LO.

\section{Appendix. Controlling the tangent step by \\ primal barrier function}\label{sc-App}
\SetEQ

In this section, we derive the lower bounds for the tangent step, using the second representation (\ref{eq-OGrow}$_2$) of the growth of the functional proximity measure. For that, we need to estimate derivatives of function $g(\cdot)$, which is a reduction of the primal barrier $F(\cdot)$ onto the set ${\cal F}_p$. For any $x \in \inter K$ and $s \in \E^*$, denote
\beq\label{def-SNorm}
\ba{rcl}
\| s \|^*_x & = & \max\limits_{h \in \E} \{ \la s, h\ra: \; A h = 0, \; \| h \|_x \leq 1 \}.
\ea
\eeq
This definition admits a closed-form representation.
\BL\label{lm-Rep}
For any $x \in \inter K$ and $s \in \E$, we have
\beq\label{eq-Rep}
\ba{rcl}
\| s \|_x^* & = & \| s - A^* u \|_x, \quad u \; = \; [A [\nabla^2 F(x)]^{-1} A^*]^{-1}A [\nabla^2 F(x)]^{-1} s.
\ea
\eeq
Therefore,
\beq\label{eq-Rep1}
\ba{rcl}
(\|s \|^*_x)^2 & = & \| s \|_x^2 - \la [A [\nabla^2 F(x)]^{-1} A^*]^{-1}A [\nabla^2 F(x)]^{-1} s, A [\nabla^2 F(x)]^{-1} s \ra.
\ea
\eeq
\EL
\proof
Denoting by $h$ the optimal solution of the problem (\ref{def-SNorm}) and by $u \in \R^m$ and $\lambda >0$ the optimal dual multipliers, we have
$$
\ba{rcl}
s & = & \lambda \nabla^2F(x) h + A^* u, \quad A h = 0, \quad \la \nabla^2 F(x) h, h \ra = 1.
\ea
$$
Hence, $h = {1 \over \lambda} [\nabla^2F(x)]^{-1}(s - A^* u)$, which gives us $\lambda^2 = \| s - A^* u \|^2_x$ and 
$$
\ba{rcl}
u & = & [A [\nabla^2 F(x)]^{-1} A^*]^{-1}A [\nabla^2 F(x)]^{-1} s.
\ea
$$
Thus, $\la s, h \ra = \la s - A^* u, h \ra = \lambda$.
\qed
Note that for any $\hat x \in K$ with $\|x - \hat x \|_x \leq \beta$, we have $\| s \|^*_{\hat x} \refLE{eq-Hess} {1 \over 1 - \beta} \| s \|^*_x$ for all $s \in \E^*$.

Now we can define the barrier parameter of function $g(\cdot)$:
$\nu_g \Def \sup\limits_{x \in \rint {\cal F}_p} \left(\| \nabla F(x) \|_x^* \right)^2$.

For method (\ref{met-DBS}), we need the following relation.
\BL\label{lem-DXK}
At the predictor step of method (\ref{met-DBS}), we have
\beq\label{eq-DXK}
\ba{rcl}
\| \Delta x_k \|_{x_k} & = & t_k \| c \|^*_{x_k},
\ea
\eeq
where $x_k = - {1 \over t_k} \nabla F_*(s_k)$.
\EL
\proof
Indeed, $\Delta x_k \refEQ{met-DBS} - \nabla^2 F_*(w^*_k) (\Delta s_k  + \hat s_k)$ and
$$
\ba{rcl}
\hat s_k  + \Delta s_k  & = & s_k - A^* (d_k+\Delta y_k) \; \refEQI{eq-FormASD}{1} \;
s_k - A^*[\nabla^2 \zeta(y_k)]^{-1}( \nabla \psi_{t_k}(y_k) + t_k b)  \\
\\
& = & s_k - A^*[\nabla^2 \zeta(y_k)]^{-1}\nabla \zeta(y_k).
\ea
$$
Since $\nabla \zeta(y_k) = - A \nabla F_*(s_k) \refEQI{eq-HomHX}{1} A \nabla^2 F_*(s_k) s_k$ and $\nabla^2 \zeta(y_k) = A \nabla^2F_*(s_k) A^*$, we get
$$
\ba{rcl}
\hat s_k  + \Delta s_k  & = & \Big(I - A^*[A \nabla^2F_*(s_k) A^*]^{-1}A \nabla^2 F_*(s_k) \Big) (c - A^* y_k)\\
\\
& = & \Big(I - A^*[A \nabla^2F_*(s_k) A^*]^{-1}A \nabla^2 F_*(s_k) \Big)c.
\ea
$$
In view of (\ref{eq-Hom12}$_2$) and (\ref{eq-FDual}), we have $\nabla^2 F_*(s_k) = t_k^2 [\nabla^2 F(x_k)]^{-1}$. Thus, since $w^*_k \refEQ{met-DBS} \sqrt{t_k} s_k$, we get the following equality:
$$
\ba{rcl}
\Delta x_k & \refEQI{eq-Hom12}{2} & - t_k [\nabla^2 F(x_k)]^{-1} \Big(I - A^*[A [\nabla^2 F(x_k)]^{-1} A^*]^{-1}A [\nabla^2 F(x_k)]^{-1} \Big)c.
\ea
$$
It remains to use representation (\ref{eq-Rep}).
\qed

Now we can justify an upper bound for the size of $\Delta x_k$ with respect to the point $\hat x_k$:
$$
\ba{rcl}
(1-\beta) \| \Delta x_k \|_{\hat x_k} & \refLE{eq-HDist} & \| \Delta x_k \|_{x_k} \; \refEQ{eq-DXK} \; t_k \| c \|^*_{x_k} \; \leq \; 
\| \nabla f_{t_k}(\hat x_k) \|^*_{x_k} + \| \nabla F(\hat x_k) \|^*_{x_k}\\
\\
& \refLE{eq-HDist} & {1 \over 1 - \beta} \Big[ \lambda_{f_{t_k}}(\hat x_k) + \lambda_g(\hat x_k) \Big] \; \stackrel{{\bf L.\ref{lm-LF}}}{\leq} \; {1 \over 1 - \beta} \Big[ \left({\beta \over 1 - \beta} \right)^2 + \sqrt{\nu_g} \Big].
\ea
$$
Thus, we have proved the following bound:
\beq\label{eq-DXKP}
\ba{rcl}
\| \Delta x_k \|_{\hat x_k} & \leq & 
{1 \over (1 - \beta)^2} \Big[ \left({\beta \over 1 - \beta} \right)^2 + \lambda_g(\hat x_k) \Big] \; \leq \; \check c_{\beta} \Def
{1 \over (1 - \beta)^2} \Big[ \left({\beta \over 1 - \beta} \right)^2 + \sqrt{\nu_g} \Big].
\ea
\eeq

Let us estimate now the growth of the functional proximity measure for self-scaled cones.
$$
\ba{rcl}
\xi_{\hat z_k}(\alpha) & \refEQI{eq-OGrow1}{2} & F(\hat x_k + \alpha \Delta x_k) + F\left( \hat x_k - {\alpha \over 1 - \alpha} \Delta x_k \right) - 2 F(\hat x_k)\\
\\
& \refLE{eq-SCF} & \la \nabla F(\hat x_k), \alpha \Delta x_k - {\alpha \over 1 - \alpha} \Delta x_k \ra + \omega_*(\alpha \check c_{\beta}) + \omega_*\left({\alpha \over 1 - \alpha} \check c_{\beta} \right)\\
\\
& \refLE{eq-GrowDY} & {\alpha^2 \over 1 - \alpha} \lambda_g(\hat x_k)  \| \Delta x_k \|_{\hat x_k} + \omega_* \left( {\alpha \check c_{\beta}(2 - \alpha) \over 1 - \alpha} \right)\\
\\
& \refLEQ{def-SCB1}{eq-DXKP} & {\alpha^2 \over 1 - \alpha} \sqrt{\nu_{g}} \check c_{\beta} + \omega_* \left( {2\alpha \check c_{\beta}\over 1 - \alpha} \right) \; \leq \; {(1-\beta)^2  \over 1 - \alpha}(\alpha \check c_{\beta})^2 + \omega_* \left( {2\alpha \check c_{\beta}\over 1 - \alpha} \right).
\ea
$$

As in the proof of Theorem \ref{th-DBS}, we conclude that the step size $\alpha_{k_1}$ at the predictor iteration $k_1$ satisfies inequality
\beq\label{eq-DProg1}
\ba{rcl}
\delta & \leq & \xi_{\hat z_{k_1}}(\alpha_{k_1})  \; \leq \; {(1-\beta)^2  \over 1 - \alpha_{k_1}}(\alpha_{k_1} \check c_{\beta})^2 + \omega_* \left( {2\alpha_{k_1} \check c_{\beta}\over 1 - \alpha_{k_1}} \right)
\ea
\eeq
Let us assume that $\alpha_{k_1} \check c_{\beta} \leq {2 \over (1 - \beta)^2}$. Then
${(1-\beta)^2  \over 1 - \alpha_{k_1}}(\alpha_{k_1} \check c_{\beta})^2 \leq {2\alpha_{k_1} \check c_{\beta}\over 1 - \alpha_{k_1}}$, and we get
$$
\ba{rcl}
\delta & \refLE{eq-DProg1} & - \ln \left(1 - {2\alpha_{k_1} \check c_{\beta}\over 1 - \alpha_{k_1}} \right) \; \leq \; {{2\alpha_{k_1} \check c_{\beta}\over 1 - \alpha_{k_1}} \over 1  - {2\alpha_{k_1} \check c_{\beta}\over 1 - \alpha_{k_1}}}.
\ea
$$
Hence, ${2\alpha_{k_1} \check c_{\beta}\over 1 - \alpha_{k_1}} \geq {\delta \over 1 + \delta}$. This means that ${\alpha_{k_1} \over 1 - \alpha_{k_1}} \geq {\delta \over 2 \check c_{\beta}(1 + \delta)}$, and therefore
$\alpha_{k_1}  \geq {\delta  \over 2 \check c_{\beta} (1+\delta) + \delta}$.
Combining this observation with our temporary assumption, we conclude that 
$$
\ba{rcl}
\alpha_{k_1}  & \geq & \min\Big\{ {2 \over (1-\beta)^2 \check c_{\beta}}, {\delta  \over 2 \check c_{\beta} (1+\delta) + \delta} \Big\} \; = \; {\delta  \over 2 \check c_{\beta} (1+\delta) + \delta} \; \geq \; {\delta  \over (2 \check c_{\beta}+1) (1+\delta)}.
\ea
$$
Hence, as in the proof of Theorem \ref{th-DBS}, we can find the following bound for the step size:
$$
\ba{rcl}
\alpha_{k_1} - {2 \over \sqrt{\nu}} \beta & \geq & {1 \over 2 \check c_{\beta}+1} \left( {\delta  \over 1+\delta} - {2 \beta (2 \check c_{\beta}+1) \over \sqrt{\nu_{g}}} \right) \; \refGE{eq-DefNu} \; {1 \over 2 \check c_{\beta}+1} \left( {\delta  \over 1+\delta} - {2 \beta (3+\beta) \over 1 - \beta} \right).
\ea
$$
Now we can present the main fact on the performance of method (\ref{met-DBS}) for symmetric cones. Its justification is the same as that one for Theorem \ref{th-DBS}. If
${\delta  \over 1+\delta} \geq {4 \beta (3+\beta) \over 1 - \beta}$, then 
$$
\ba{rcl}
\alpha_{k_1} - {2 \over \sqrt{\nu}} \beta & \geq & {\delta \over 4 ( \check c_{\beta}+1)(1+\delta)},
\ea
$$
and we get that the penalty coefficient $t_k$ grows at least as 
$\exp \Big\{ {k \over O( \sqrt{\nu_g})} \Big\}$. 
\qed

\end{document}